\documentclass[12pt]{article}
\usepackage{amsthm,amssymb,hyperref}
\input xypic \xyoption{all}
\title{Asymptotic behaviour and \\ the moduli space of doubly-periodic instantons}
\author{Olivier Biquard \\ \'Ecole Polytechnique \\ \textsc{cmat, umr} 7640 du
\textsc{cnrs} \\ 91128 Palaiseau Cedex, France \\ \\ and \\ \\ Marcos Jardim \\ 
Yale University \\ Department of Mathematics \\ 10 Hillhouse Avenue \\ 
New Haven, CT 06520-8283 USA}

\newcommand{\pf}{{\em Proof. }} \newcommand{\pfend}{\qed}
\newcommand{\seta}{\rightarrow} \newcommand{\as}{\pm\ksi_0}
\newcommand{\torus}{T\times\cpx} \newcommand{\tproj}{T\times\proj}
\newcommand{\del}{\overline{\partial}} \newcommand{\ksi}{\xi}
\newcommand{\bR}{\mathbb{R}} \newcommand{\real}{\mathbb{R}} 
\newcommand{\cpx}{\mathbb{C}} \newcommand{\zed}{\mathbb{Z}}
\newcommand{\proj}{\mathbb{P}^1} \newcommand{\hh}{\mathbb{H}}
\newcommand{\ee}{\mathcal{E}} \newcommand{\dual}{\hat{T}} 
\newcommand{\ff}{\mathcal{F}} \newcommand{\vv}{\mathcal{V}} 
\newcommand{\cM}{\mathcal{M}} \newcommand{\hmod}{\hat{\mathcal{M}}}
 \newcommand{\cF}{\mathcal{F}}  \newcommand{\oo}{\mathcal{O}}
\newcommand{\cG}{\mathcal{G}} \newcommand{\cA}{\mathcal{A}} 
\newcommand{\hI}{\hat{I}} \newcommand{\su}{\otimes\mathfrak{su}}
\newcommand{\osl}{\otimes\mathfrak{sl}} \newcommand{\ou}{\otimes\mathfrak{u}}

\newtheorem{thm}{Theorem}[section] \newtheorem{lem}[thm]{Lemma}
\newtheorem{prop}[thm]{Proposition} 

\theoremstyle{remark}
\newtheorem{rem}[thm]{Remark}
\newtheorem*{claim}{Claim}

\begin{document}

\maketitle

\begin{abstract}
We study  doubly-periodic instantons, i.e. instantons on the product
of a 1-dimensional complex torus $T$ with a complex line $\mathbb{C}$,
with quadratic curvature decay.
We determine the asymptotic behaviour of these instantons, 
constructing new asymptotic invariants. We show that the underlying holomorphic 
bundle extends to $T\times\mathbb{P}^1$.
The converse statement is also true, namely a holomorphic bundle on
$T\times\mathbb{P}^1$ which is flat on the torus at infinity, and 
satisfies a stability condition, comes from a doubly-periodic instanton. 
Finally, we study the hyperk\"ahler geometry of the moduli space of doubly-periodic 
instantons, and prove that the Nahm transform previously defined by the second 
author is a hyperk\"ahler isometry with the moduli space of certain meromorphic 
Higgs bundles on the dual torus.
\end{abstract}

\baselineskip18pt \newpage
\tableofcontents \newpage

\section*{Introduction and statement of the results}
The aim of this paper is to understand the analytical properties
of certain finite energy solutions of the Yang-Mills
anti-self-dual equations over $\torus$. These so-called {\em
extensible doubly-periodic instantons} have been studied by the
second author in \cite{J1,J2,J3}, where they were shown to be
equivalent to certain singular solutions of Hitchin's equations
over an elliptic curve via a construction known as the {\em Nahm
transform}. The present paper grew from questions raised in the
works mentioned above.

More precisely, consider an $SU_2$ bundle $E\seta\torus$.
The instanton connections $A$ considered in \cite{J1,J2} satisfied the
following hypothesis:
\begin{enumerate}
\item \emph{quadratic curvature decay}: $|F_A|=O(r^{-2})$ with
respect to the Euclidean metric on $\torus$;
\item \emph{extensibility}: there is a holomorphic rank two vector bundle
$\ee\seta\tproj$ with trivial determinant such that
$\ee|_{T\times(\proj\setminus\{\infty\})}\simeq(E,\del_A)$, where
$\del_A$ is the holomorphic structure on $E$ induced by $A$;
\end{enumerate}
where $w$ is a coordinate in the complex line, and by the notation
$O(|w|^\gamma)$ we mean the set of functions on $\cpx$ such that:
$\lim_{|w|\seta\infty}|f(w)|/|w|^\gamma < \infty$.

One of the goals of this paper is to prove that the technical hypothesis
of extensibility is actually a consequence of the anti-self-duality
equation, and more generally to understand completely the behaviour at
infinity of all instantons with quadratic curvature decay.

\smallskip\noindent\textbf{Model solutions.}
Special solutions of the anti-self-duality equations may be obtained
by restricting to torus invariant connections. Such instantons
come from solutions $(B,\psi)$ of Hitchin's equations on $\cpx$
$$ \left\{ \begin{array}{l}
F_B+[\psi,\psi^*]=0 \\
\del_B\psi=0
\end{array} \right. $$
in the following way. Recall that $B$ is a $SU_2$-connection on
$\cpx$, and $\psi$ is a (1,0)-form with values in $\mathfrak{sl}_2$.
Let $\psi=\frac{1}{2}(\psi_0+i\psi_1)dw$, and consider the connection
(where $x$ and $y$ are coordinates on $T$):
$$ A_0 = B + \psi_0 dx + \psi_1 dy $$
which is a torus invariant instanton. Assuming that $|F_{A_0}|=O(r^{-2})$,
the asymptotic behavior of solutions $(B,\psi)$ is given by one of the
following models:
\begin{eqnarray}
B = d \qquad & \quad & \psi=
\left(\begin{array}{cc} \lambda & 0 \\ 0 & -\lambda
      \end{array}\right)dw \label{ex1} \\
B=d+\left(\begin{array}{cc} \alpha & 0 \\ 0 & -\alpha
        \end{array}\right) d\theta & \quad &
\psi=\left(\begin{array}{cc} \mu & 0 \\ 0 & -\mu
           \end{array}\right)\frac{dw}{w} \label{ex2} \\
B=d+\left( \begin{array}{cc} -1 & 0 \\ 0 & 1
         \end{array}\right)\frac{d\theta}{\ln r^2} & \quad &
\psi=\left(\begin{array}{cc} 0 & 1\\ 0 & 0
           \end{array}\right)\frac{dw}{w\ln r^2} \label{ex3}
\end{eqnarray}
where $\lambda,\mu\in\cpx$ and $-\frac{1}{2}\leq\alpha<\frac{1}{2}$.
The solutions of examples (\ref{ex1}) \& (\ref{ex2}) can be
superimposed, and such superpositions
are called the {\em semisimple} solutions. On the other hand,
solutions of example (\ref{ex3}) cannot be superimposed with the others;
these are called the {\em nilpotent} solutions, and can only exist
when  $\lambda=\mu=\alpha=0$.
The torus invariant instanton is then given by, in the semisimple
case:
$$ A_0=d+i\left(\begin{array}{cc}a_0 & 0 \\ 0 & - a_0
\end{array}\right) $$
with
$$ a_0 = \lambda_1 dx + \lambda_2 dy + (\mu_1 \cos\theta - \mu_2 \sin\theta)
\frac{dx}{r} + (\mu_1 \sin\theta + \mu_2\cos\theta) \frac{dy}{r} + \alpha d\theta ; $$
while in the nilpotent case, we have:
$$ A_0 = d+i\left(\begin{array}{cc} -1 & 0 \\ 0 & 1 \end{array}\right)
   \frac{d\theta}{\ln r^2} + \frac{1}{r\ln r^2}
   \left(\begin{array}{cc} 0 & e^{-i\theta}(dx-idy) \\ -e^{i\theta}(dx+idy) & 0
   \end{array}\right) $$
and note that the curvature is $O(r^{-2}|\ln r^2|^{-2})$.

Remark that the connection $A_0$ has a flat limit over the torus at infinity,
$$
d+i\left(\begin{array}{cc}\lambda_1 dx+\lambda_2 dy & 0 \\ 0 &
-\lambda_1 dx - \lambda_2 dy \end{array}\right) ,
$$
and one can prove that such flat limit for a connection $A$ exists as
soon as $|F_A|=O(r^{-1-\epsilon})$; the flat limit underlies a
holomorphic vector bundle $L_{\xi_0}\oplus L_{-\xi_0}$, where the
elements of the dual torus $\pm\xi_0\in\dual$ are called the
\emph{asymptotic states} of the connection.

We show that the three standard examples above completely describe
the behavior at infinity of doubly-periodic instantons with
quadratic curvature decay:
\begin{thm} \label{asymp.par}
Let $A$ be a doubly-periodic instanton with curvature $O(r^{-2})$.
Then there is a gauge near infinity such that
$$ A=A_0+a , $$
where $A_0$ is one of the previous models, and, for some
$\delta>0$, in the semisimple case:
$$ |a|=O\bigg(\frac{1}{r^{1+\delta}}\bigg) , \quad
   |\nabla_{A_0}a|=O\bigg(\frac{1}{r^{2+\delta}}\bigg) ; $$
in the nilpotent case:
$$ |a|=O\bigg(\frac{1}{r (\ln r)^{1+\delta}}\bigg) , \quad
   |\nabla_{A_0}a|=O\bigg(\frac{1}{r^2 (\ln r)^{2+\delta}}\bigg) . $$
\end{thm}

In the case where the limit at infinity of $A$ is non
trivial, one can prove the theorem under the weaker assumption that
the curvature is $O(r^{-1-\epsilon})$; this condition is very close
to the finite energy condition, and it is natural to suppose that the
theorem actually describes the behaviour of all finite energy
instantons. The instantons we will use (for example, those coming
from the inverse Nahm transform) have quadratic curvature decay,
so that this hypothesis is sufficient for our applications.

The theorem, to be proved in section \ref{ab}, provides a complete
characterization of the instanton parameters which are invariant
under $L^2$ deformations. The parameter $\lambda$ is equivalent to
the asymptotic states $\as$. The two remaining parameters are new:
$\alpha$ is called the {\em limiting holonomy} of the instanton $A$,
while $\mu$ is called the {\em residue}. The motivation for the
latter nomenclature will be made clear latter on. Notice that, in
contrast with the instanton number (see below) and the asymptotic states, the
limiting holonomy and the residues are defined only for
anti-self-dual connections.

\smallskip\noindent\textbf{Instantons and holomorphic bundles.}
We are now ready to state our second main result, which in particular
solves the extensibility problem.
Recall \cite{J2} that the {\em instanton number} $k$ of the
doubly-periodic instanton $A$ is defined by the formula:
$$ k = \frac{1}{8\pi^2} \int_{\torus} |F_A|^2 $$
as usual.

\begin{thm} \label{extn}
There is a 1-1 correspondence between the following objects:
\begin{itemize}
\item $SU_2$-doubly-periodic instanton connections with quadratic
curvature decay and fixed asymptotic parameters $(k,\as,\alpha)$;
\item $\alpha$-stable, rank two holomorphic vector bundles $\ee\seta\tproj$ with
trivial determinant such that $c_2(\ee)=k$ and
$\ee|_{T\times\{\infty\}}=L_{\xi_0}\oplus L_{-\xi_0}$.
\end{itemize} \end{thm}

The stability condition of the statement is a variant of the
stability condition for \emph{parabolic bundles};
the degree is calculated with respect to a non ample class
(the fundamental class of the torus). The precise definition will be
given in section \ref{holo.extn}, where this result is proved.

In a broader context, theorem \ref{extn} can be seen as the
analog of Donaldson's correspondence between instantons on
$\real^4$ and framed holomorphic bundles over $\mathbb{P}^2$
\cite{D,Bu}. In this last case, no stability condition is needed in
order to produce an instanton, while in the case of a compact surface,
stability (with respect to an ample class) is necessary. Thus, in some
sense, our stability criterion goes midway between these two situations.

\smallskip\noindent\textbf{Moduli space.}
We then pass to the analytical construction of the moduli
space of doubly-periodic instantons. We prove:

\begin{thm} \label{modthm}
The moduli space of doubly-periodic instantons with fixed
instanton number $k$ and asymptotic parameters $(\as,\alpha,\mu)$
is a smooth hyperk\"ahler manifold of real dimension $8k-4$.
\end{thm}

Of course, this theorem is interesting only if the moduli space is not empty.
Fortunately, as mentioned in \cite{J1}, existence of doubly-periodic
instantons for generic values of the parameters $(k,\as,\alpha,\mu)$
is guaranteed via the Nahm transform (see below) of meromorphic Higgs
bundles over $\dual$, whose existence follows from Simpson \cite{S}
among others; theorem \ref{asymp.par} puts these instantons in our
moduli spaces. Another equivalent, probably more direct, way for guaranteeing
existence is of course theorem \ref{extn}.
See also section \ref{holo.extn} for some cases where the moduli space is
empty, and section \ref{sec.moduli} for a description of the $k=1$
moduli space.

\smallskip\noindent\textbf{Nahm transform.}
Finally we revisit the Nahm transform of doubly-periodic instantons
defined in \cite{J2} with two main objectives in mind. Before
explaining what these objectives are, let us say a few words about
the Nahm transform.

Here we restrict to the semisimple case, since Nahm transform was
defined only in this case.
Recall from \cite{J2} (see also part \ref{part3})
the Nahm Transform is a 1-1 correspondence
between irreducible, doubly-periodic instantons and
certain meromorphic Higgs pairs $(B,\Phi)$ on a bundle $V$ over the dual torus
$\dual$. The rank of $V$ is given by the instanton number.
The Higgs field $\Phi$ has simple poles at the two points
corresponding to the asymptotic states $\as$. Moreover, $\Phi$
has semisimple residues of rank one if $\xi_0\neq-\xi_0$, and
rank two otherwise. We denote by ${\rm Res}\Phi(\as)$ the residue
of the Higgs field at the singular point $\as$.

Thus, it is natural to ask how are the new asymptotic parameters
defined by theorem \ref{asymp.par} interpreted in terms of the Nahm
transformed Higgs pair. This question in answered in section
\ref{ap.nt}, and the precise statement is given in theorem
\ref{ap.nt.thm}. As expected from the general principle
{\em Nahm transform is a non-linear Fourier transform}, the
asymptotic behavior is converted into singularity behavior.

It is well known that the moduli space of Higgs pairs on a Riemann
surface is hyperk\"ahler \cite{H}; for the moduli space of Higgs pairs
with fixed singularities at the punctures, this follows from \cite{B3}.
The second goal can now be summarized in our last result:
\begin{thm} \label{iso.thm}
The Nahm transform of doubly-periodic instantons is a
hyperk\"ahler isometry.
\end{thm}
Note that similar results have been proved for the other
well-known examples of Nahm transform: the ADHM construction,
see \cite{M}; the duality between monopoles and solutions of
Nahm equations, see \cite{N}; and the Fourier-Mukai transform
of instantons over 4-tori, see \cite{BVB}. Indeed, it is
reasonable to expect that such result holds for {\em any} Nahm
transform.

\smallskip\noindent\textbf{Outline.}
The paper is divided in three parts. 
The first part is technical: we study the
asymptotic behavior of connections on $E$ with quadratic curvature decay,
but which are not necessarily anti-self-dual;
the technical goal is the construction of a {\em partial
Coulomb gauge} (theorem \ref{coulomb.gauge}).
In the second part, we obtain theorems \ref{asymp.par}, \ref{extn} and \ref{modthm}.
Finally, the third part deals with the Nahm transform aspects of the paper.

\smallskip\noindent\textbf{Acknowledgements.}
The second author would like to thank the \'Ecole Polytechnique for
its support, and Antony Maciocia for useful conversations.

\part{Connections with quadratic curvature decay} \label{part1}
In this part, we study the behaviour at infinity of (not necessarily
anti-self-dual) connections with quadratic curvature decay on a
$SU_2$-bundle $E$ on $T\times\bR^2$. Such connections will have a limit flat
connection $\Gamma$ on the torus at infinity $T_\infty$, which
decomposes $E|_{T_\infty}$ as a sum of two flat line bundles
$L_{\xi_0}\oplus L_{-\xi_0}$; when $L_{\xi_0}^2=0$, we can reduce to
the case where $L_{\xi_0}=0$ by globally tensoring $E$ with
$L_{\xi_0}$; therefore we will always suppose that
\begin{equation}\label{Lxi2}
\textnormal{either }L_{\xi_0}^2\neq 0, \quad\textnormal{ or } L_{\xi_0}=0 .
\end{equation}

Over any torus $T$, we consider the $L^2$-orthogonal decomposition
\begin{equation}\label{decomposition01}
End(E) = (\ker \nabla_\Gamma) \oplus (\ker \nabla_\Gamma)^\perp
\end{equation}
and we decompose accordingly any section $u$ of $End(E)$ as
\begin{equation}
u = u_\Gamma + u_\perp .
\end{equation}
If we write explicitly $\Gamma=d+\gamma$, with
$$
\gamma=\left(\begin{array}{cc}\lambda_1 & 0 \\ 0 & -\lambda_1
             \end{array}\right) dx +
       \left(\begin{array}{cc}\lambda_2 & 0 \\ 0 & -\lambda_2
             \end{array}\right) dy,
$$
then, in view of (\ref{Lxi2}),
$\ker\nabla_\Gamma$ is described as the $T$-invariant sections of
$\ker\gamma$; if $\gamma$ is nontrivial, these are reduced to
$T$-invariant diagonal matrices.

The technical goal of this part is a {\em partial Coulomb
gauge} on the $a_\perp$ part of a connection $A=\Gamma+a$
with curvature $O(r^{-2})$. More precisely, let $V_R\subset\cpx$
denote the complement of a disc of radius $R$ centered at the origin.

\begin{thm} \label{coulomb.gauge}
Given a constant $\eta>0$, there exists $R$ sufficiently large
such that if $A$ is a doubly-periodic connection satisfying
$\sup_{r\geq R} \big( r^2 |F_A| \big) \leq \eta$, then there is a
gauge $g$ on $T\times V_R$ such that
$g(A)=\Gamma+a_\Gamma+a_\perp$, with:
\begin{eqnarray*}
\textnormal{(i)} & & d^*_{\Gamma+a_\Gamma} a_\perp = 0 ; \\
\textnormal{(ii)} & & \partial_r \lrcorner a_\perp (r=R) = 0 ; \\
\textnormal{(iii)} & & \| r^2 F_{\Gamma+a_\Gamma} \|_{C^0}
 + \| r^{2-\epsilon} a_\perp \|_{C^0} \leq C \cdot \| r^2 F_A \|_{C^0}  .
\end{eqnarray*} \end{thm}
Note that gauge transformations $g=g_\Gamma$ preserve the Coulomb
gauge constructed in this theorem. This kind of partial gauge
fixing reminds of R{\aa}de's fibered Hodge gauge \cite{R}.

\begin{rem}\label{better-decay}
Actually, if $\Gamma$ is nontrivial, the proof gives a Coulomb
gauge under a weaker bound on the curvature, namely
$|F|=O(r^{-(1+\epsilon)})$; this condition is very close to the finite
energy condition, since $r^{-\delta}$ is in $L^2$ when $\delta>1$.
\end{rem}

\section{Limit flat connection}

Our first task is to establish the existence of a flat limit
connection $\Gamma$ for every connection $A$ with quadratic
curvature decay:

\begin{prop}\label{flat-limit}
Suppose that the connection $A$ on $T\times\bR^2$ satisfies $$
|F_A| \leq \frac{c_1}{r^2}. $$ Then $A$ has a flat limit $\Gamma$
on $T$ at infinity, and there exists a sequence of connections
$A_j$, such that
\begin{enumerate}
\item $|F_{A_j}| \leq c_2/r^2$ ;
\item $A_j$ is gauge equivalent to $A$ on $\{r\leq j\}$ ;
\item $A_j=\Gamma$ on $\{r\geq 2j\}$.
\end{enumerate}
\end{prop}

\begin{rem} This proposition remains true if the curvature is
$O(r^{-(1+\epsilon)})$. \end{rem}

\begin{proof}
We begin by proving the existence of the flat limit $\Gamma$. Take
a radial gauge $$ A = d + a_\theta d\theta + a_x dx + a_y dy $$
for $A$; from the bound on the curvature, we deduce
\begin{equation}\label{radial-bounds}
|\partial_r a_x|+|\partial_r a_y| = O(r^{-2}) , \quad |\partial_r
a_\theta| = O(r^{-1}) ;
\end{equation}
from this we deduce that $a_x$ and $a_y$ have limits
$a^{\infty}_x(\theta,x,y)$ and $a^\infty_y(\theta,x,y)$ when $r$
goes to infinity; moreover, the bound on the curvature implies
that for each $\theta$, the connection
$d+a^{\infty}_x(\theta)dx+a^\infty_y(\theta)$ is flat on $T$. It
remains to see that it is independent of $\theta$: for this we
pick a base point in $T=S^1\times S^1$ and prove that the
monodromies along the two circles remain conjugate when $\theta$
varies; this is a consequence of the bound on the curvature and
the following lemma (see for example \cite[lemma 1]{B1}):

\begin{lem}\label{monodromy-curvature}
Suppose we have a connection $A$ on $[0,1]\times S^1$, and $m(t)$
is the monodromy of $A$ along the circle $\{t\}\times S^1$; note
$h(t)$ the parallel transport from the point $(0,0)$ to the point
$(t,0)$; then $$ \big|\partial_t \big(h(t)^{-1} m(t)
h(t)\big)\big| \leq \int_{\{t\}\times S^1} |F_A| .$$
\end{lem}

Therefore we have constructed a flat limit $\Gamma$ on $T$ for the
connection $A$. Now pass to the approximation statement.

\begin{claim}
On $\{r\}\times S^1 \times T$, there exists a gauge so that
$A=\Gamma+a$, $|a|\leq c/r$.
\end{claim}

This statement (a $C^0$ gauge only), can be proven by elementary
means and is left to the reader. Now, we extend radially this
gauge on $\{j\}\times S^1 \times T$ to $[j,2j]\times S^1 \times
T$, and the bounds (\ref{radial-bounds}) imply that $A=\Gamma+a$
with still $|a|\leq c/r$ on $[j,2j]$; then we choose a cutoff
function $\chi=\chi(r)$ so that $$ \chi(r\leq j)=1, \quad
\chi(r\geq 2j)=0, \quad |\partial_r\chi|\leq 2/j,$$ and define a
connection $A_j$ by
\begin{eqnarray*}
r\leq j, & \quad & A_j = A , \\ r\geq j, & \quad & A_j = \Gamma +
\chi a  ;
\end{eqnarray*}
on $r\geq j$, the curvature of $A_j$ is $$ F_{A_j} = \chi F_A +
d\chi \wedge a + (\chi^2-\chi) a\wedge a $$ and this remains
bounded by $c/r^2$ on $[j,2j]$, which means that $|F_{A_j}|$ is
uniformly bounded by $c/r^2$.
\end{proof}

\begin{rem} Actually, it is not difficult to go a bit further and to
prove that there is a global gauge in which $A=\Gamma+a$ and
$|a|=O(\ln r/r)$; this gives a result used without proof in
\cite{J2}. Of course, the result will also be a consequence of
theorem \ref{coulomb.gauge}.
\end{rem}

In the case of a torus invariant connection, we need a stronger
statement.
\begin{prop}\label{invariant-control}
Under the hypotheses of proposition \ref{flat-limit}, if $A=d+a$
with $a=a_\Gamma$ (in particular $A$ is torus invariant), then
there is a gauge such that $$A=\Gamma+a+b , $$ where $d+a$ is a
connection on $\bR^2$ and $b=b_x dx+b_y dy$ a 1-form along $T$,
satisfying $$ |b|\leq\frac{c_3}{r}, \quad
|\nabla_{\Gamma+a}b|\leq\frac{c_2}{r^2},$$ and $$ a=i
\left(\begin{array}{cc}\alpha(r) & 0 \\ 0 &
-\alpha(r)\end{array}\right)
     d\theta + b $$
with $$ |\partial_r\alpha|+|b|\leq c_3/r ,
   \quad \sup_j j^{2-2/p} \|\nabla b\|_{L^p(j\leq r\leq 2j)} \leq c_3 . $$
\end{prop}
The meaning is that we want a gauge with not only a $C^0$ bound,
but also a $C^1$ bound; actually this is not possible (because
elliptic regularity does not hold in $C^k$ spaces) and this
explains why we use $L^p$ derivatives instead. So the proposition
must be considered as a regularization of the connection. The
standard way to obtain this is to use Hodge gauges in order for
the curvature to become an elliptic equation: locally Uhlenbeck's
theorem provides the required statement, but the glueing is not
easy, especially on a non simply connected manifold. We present
here a proof based on the following lemma, which is a consequence
of the Hodge gauge constructed in \cite[theorem 1]{B1}:
\begin{lem}\label{hodge-cylinder}
Any connection $A$ on $[0,1]\times S^1$, with $\|F_A\|_{L^p}$
sufficiently small, is gauge equivalent to a connection $d+i\alpha
d\theta+a$ with $\|a\|_{L^{1,p}} \leq c \|F_A\|_{L^p}$, where
$\alpha$ is a diagonal matrix, with coefficients in $[0,1[$, such
that $\exp(-2\pi i \alpha)$ is the monodromy of $A$ along the
circle $\{0\}\times S^1$. \qed
\end{lem}

\begin{proof}[Proof of proposition \ref{invariant-control}]
If $A$ is torus invariant, then a torus invariant gauge
transformation $g$ acts on $b$ only by $gbg^{-1}$, and the bounds
on the curvature immediately imply the required bounds on $b$.
Therefore, we are reduced to look at a connection $d+a$ on
$\bR^2$.

Now note that the region $r\geq 1$ is conformally equivalent to the
half-cylinder {$\bR_+\times S^1$} (with coordinate $t=\ln r$); in
the rest of the proof we will use only the flat metric on the
cylinder. The bound on the curvature becomes $|F_A|\leq c_2$;
eventually pulling back $A$ using the transformation $t
\rightarrow \lambda t$ with $\lambda$ sufficiently small, we may
suppose that $c_2$ is very small. This means that we are now able
to use lemma \ref{hodge-cylinder}, for some $p$ very big, to
produce on each $[j-1,j+1]\times S^1$ a gauge $g_j$ so that $$
g_j(A) = d_{\alpha_j} + a_j , \quad
  \|a_j\|_{L^p} + \|(\nabla+i\alpha_j d\theta)a_j\|_{L^p} \leq c c_2 . $$
We perform recursively diagonal gauge transformations with
coefficients of type $\exp(ik\theta)$ ($k$ integer) so that we
have $$ |\alpha_{j+1}-\alpha_j| < c_2 ; $$ this is possible
because of lemma \ref{monodromy-curvature}, and the operation does
not affect the bound on $a_j$ (but we have only $ |\alpha_j| \leq
c_2 j $). We want to glue together these local gauges: the
transition $h_j=g_{j+1}\cdot g_j^{-1}$ satisfies $$ dh_j +
[\alpha_j,h_j] = h_j\cdot a_j -
(\alpha_{j+1}-\alpha_j+a_{j+1})\cdot h_j ; $$ the RHS is controled
by $c c_2$, and this implies that $h_j$ is very close to some
$\tilde{h}_j\cdot(\theta)$ in the kernel of $d+\alpha_j$;
replacing $g_{j+1}$ by $\tilde{h}_j\cdot g_{j+1}$, we now may
suppose that the transition $g_{j+1}\cdot g_j^{-1}$ is close to
the identity (in $L^{2,p}$ norm), and a standard argument now
enables us to glue together all these gauges: for a similar
argument, see \cite[pages 447--8]{B1}. If we choose diagonal
matrices $\alpha(t)$ so that $$ \alpha(j)=\alpha_j, \quad
|\partial_t\alpha|\leq 2c_2,$$ we finally get a gauge
$d+i\alpha(t)d\theta+b$, with $$ \|b\|_{L^p([j-1,j+1])} +
\|(\nabla+i\alpha(t)d\theta)b\|_{L^p([j-1,j+1])} \leq c c_2 . $$
Sobolev embedding implies that $\|b\|_{C^0}$ is controled as well;
translating back these bounds in the metric of $\bR^2$, we get the
proposition.
\end{proof}

\begin{rem}
The proof of proposition \ref{invariant-control} becomes certainly
easier if $A$ is abelian (which is the case if the limit $\Gamma$
is regular), since in this case, it is easy to produce a global
Hodge gauge.
\end{rem}

\begin{rem}
In general, we are unable to prove proposition \ref{invariant-control}
if the curvature is only $O(r^{-(1+\epsilon)})$: this is because, in
order to use get a controled gauge on $\bR_+\times S^1$, we need the
curvature to be bounded;
if $\Gamma$ is nontrivial, the problem becomes abelian, and then it is
easy to construct a global Hodge gauge on $\bR_+\times S^1$, from
which the proposition follows easily (and one gets a bound in
$O(r^{-\epsilon})$ on $a$).
\end{rem}

%-----------------------------------------------------------------------------

\section{The linear problem} \label{estimates}

In this section we study the linear analysis on the $(\ker
\nabla_\Gamma)_\perp$ part for the Laplacian operator
$d_\Gamma^*d_\Gamma$ acting on 0-forms and the deformation
operator $d^+_\Gamma+d^*_\Gamma$ acting on 1-forms, with fixed
boundary conditions.

For this analysis, we will use the Sobolev spaces $L^{p,k}$ of
functions with $k$ derivatives in $L^p$; the weighted Sobolev
spaces $L^{p,k}_\delta$ of functions $f$ such that
$(1+r^2)^{\delta/2}f \in L^{p,k}$.

The basis of the analysis is the following simple lemma, which is
an immediate consequence of the decomposition
(\ref{decomposition01}).
\begin{lem} \label{l2}
There is a constant $c$, depending on $p$, such that on each torus
T, for any section $u$ of $End(E)$, we have:
\begin{equation}\label{est0}
\int_T | \nabla_\Gamma u_\perp |^p \geq c \int_T |u_\perp|^p .
\end{equation} \end{lem}

\vskip12pt \centerline{\bf Analysis on 0-forms} \vskip10pt

\begin{lem}\label{neumann-L2}
The Neumann problem on sections of $End(E)_\perp$ on $r\geq R$,
\begin{equation}\label{neumann}
\left\{ \begin{array}{l} \Delta_\Gamma u = v \\ \partial_r
u(r=R)=0
\end{array} \right. \end{equation}
is an isomorphism $L^{2,2}_\delta \rightarrow L^2_\delta$.
\end{lem}

\pf The solution $u$ of the Neumann problem is obtained by
minimizing the functional $$ \int \frac{1}{2} |\nabla_\Gamma u|^2
- \langle u,v \rangle $$ in the space $L^{1,2}$; the minimization
is possible because of the estimate (\ref{est0}); local elliptic
regularity gives that the $L^{1,2}$-solution actually lives in
$L^{2,2}$, and this gives the statement when there is no weight.

In the case we have a weight $\delta$, the following estimate
holds:
\begin{eqnarray*}
\int \langle \Delta_\Gamma u,u \rangle r^{2\delta} &=& \int
|\nabla_\Gamma u|^2 r^{2\delta} +
    2 \frac{\delta}{r} \langle \nabla_{\partial_r}u,u \rangle r^{2\delta}\\
&\geq& \int (1-\frac{\delta}{r}) |\nabla_\Gamma u|^2 r^{2\delta} -
       \frac{\delta}{r} |u|^2 r^{2\delta}
\end{eqnarray*}
and using (\ref{est0}) we get, if $R$ is large enough,
\begin{eqnarray*}
\|\Delta_\Gamma u\|_{L^2_\delta} \|u\|_{L^2_\delta} &\geq& \int
\langle \Delta_\Gamma u,u \rangle r^{2\delta} \\ &\geq& C
\|u\|_{L^2_\delta}^2
\end{eqnarray*}
and therefore $$ C \|u\|_{L^2_\delta} \leq \|\Delta_\Gamma
u\|_{L^2_\delta} $$ which proves that the isomorphism persists
between weighted $L^2$-spaces, at least if $R$ is large enough.

This would be enough for our applications, but one can prove
easily that the statement remains true for any $R$: because
$\Delta_\Gamma$ is an isomorphism for $R$ big enough, it remains a
Fredholm operator for any $R$ (just glue the inverse near infinity
with a parametrix on the compact part); the index is locally
constant and therefore does not depend on the weight $\delta$;
this means that it is equal to the $L^2$-index, that is $0$; now,
because the $L^2$-kernel is zero, the $L^2_\delta$-kernel is zero
if $\delta>0$; for general $\delta$, the kernel is the
$L^2_\delta$-kernel, while the cokernel is the
$L^2_{-\delta}$-kernel: as at least one of them is trivial and the
index is $0$, both are trivial. \pfend

We now want to deduce the same result in $L^p$ spaces. We need an
estimate on the solution of problem (\ref{neumann}) when $v$ is
$L^p$. After a conformal change in the Euclidean metric $g_E$, we
can pass to the cusp metric ($r=e^t$): $$
g_C=dt^2+d\theta^2+e^{2t}(dx^2+dy^2)=\frac{1}{r^2}g_E $$ The
operator $\Delta_\Gamma$ now has singular coefficients, but is
basically of the type studied in \cite{B3}, where H\"older and
$L^p$ estimates are deduced from the $L^2$-estimates. Here, the
same techniques lead to the desired result:
\begin{lem}\label{neumann-Lp}
The Neumann problem (\ref{neumann}) for 0-forms on $r\geq R$ is an
isomorphism $L^{2,p}_\delta \rightarrow L^p_\delta$ for all
weights $\delta$.
\end{lem}

\pf For the convenience of the reader, we give here a sketch of
proof for the statement, inspired from \cite[section 6]{B3}, but
written with respect to the Euclidean metric. The proof below
works for $p>2$ (the case we will use), but the statement remains
true for general $p$.

The first step is to give an elliptic estimate
\begin{equation}\label{est1}
\|u\|_{L^{2,p}_\delta[r,2r]} \leq c \big( \|\Delta_\Gamma
u\|_{L^p_\delta([\frac{1}{2}r,3r])}
      + \|u\|_{L^2_{\delta-1+2/p}([\frac{1}{2}r,3r])} \big).
\end{equation}
The weight
\begin{equation}\label{delta2}
\delta_2=\delta-1+2/p
\end{equation}
chosen for the $L^2$ space corresponds to functions with the same
order of decreasing in $r^{-\delta-2/p}$ as in the weighted $L^p$
space, but actually the proof below will give more. In order to
prove this, we remark that $$ \nabla_\Gamma = e^{-i(ax+by)} \circ
\nabla \circ e^{i(ax+by)} $$ so that if we consider $x$ and $y$ as
coordinates on $\real^2$, the equation $\Delta_\Gamma u=v$ becomes
equivalent to $$ \Delta u'=e^{i(ax+by)} v , \quad u'=e^{i(ax+by)}
u. $$ In the domain $[1,2]\times S^1\times [-1,1]^2 \subset
\real^2 \times \real^2$, we have an elliptic estimate $$
\|u'\|_{L^{2,p}} \leq c \big( \|u'\|_{L^2} + \|\Delta u'\|_{L^p}
\big) $$ which implies on the homothetic domain $[R,2R]\times
S^1\times [-R,R]^2 \subset \bR^2 \times \bR^2$ $$
R^{2-4/p}\|\nabla^2 u'\|_{L^p} \leq c \big( R^{-2} \|u'\|_{L^2} +
R^{2-4/p} \|\Delta u'\|_{L^p} \big) $$ and therefore on
$[R,2R]\times S^1\times T$ $$ R^{2-2/p}\|\nabla_\Gamma^2 u\|_{L^p}
\leq c \big( R^{-1} \|u\|_{L^2} + R^{2-2/p} \|\Delta_\Gamma
u\|_{L^p} \big) $$ which we can rewrite, still on $[R,2R]\times
S^1\times T$, $$ \|\nabla_\Gamma^2 u\|_{L^p_\delta} \leq c \big(
\|u\|_{L^2_{\delta-3+2/p}} + \|\Delta_\Gamma u\|_{L^p_\delta}
\big) $$ now the estimate (\ref{est0}) implies $$
\|\nabla^k_\Gamma u\|_{L^p} \geq c \|u\|_{L^p} ; $$ this, with
local elliptic regularity, gives the estimate (\ref{est1}).

The second step now consists in going from the $L^2$-estimates
with weights to the $L^p$-estimate. Basically, one can do the
following: let $P$ be the inverse obtained by the
$L^2$-resolution; decompose
\begin{equation}\label{vi}
v = \sum v_i
\end{equation}
where $v_i$ has support in $\exp(i/2) < r < \exp(3i)$; by the
$L^2$-resolution for the weight $\delta_2$ defined by
(\ref{delta2}), one has $$ \|Pv_i\|_{L^2_{\delta_2}} \leq c
\|v_i\|_{L^2_{\delta_2}}
   \leq c \|v_i\|_{L^p_\delta} ; $$
on the other hand, we decompose similarly $u_i=Pv_i$ as $$ u_i =
\sum_j u_{ij}, $$ and we note that the $L^2$ resolution gives the
estimate
\begin{eqnarray*}
\|u_{ij}\|_{L^2_{\delta_2}} &\leq& c e^{-\epsilon i}
\|u_{ij}\|_{L^2_{\delta_2+\epsilon}} \\ &\leq& c e^{-\epsilon i}
\|v_i\|_{L^2_{\delta_2+\epsilon}} \\ &\leq& c e^{-\epsilon (i-j)}
\|v_i\|_{L^2_{\delta_2}} ;
\end{eqnarray*}
if we choose $\epsilon$ to be $\pm \epsilon$ according to the sign
of $i-j$, we get the estimate
\begin{eqnarray*}
\|u_{ij}\|_{L^2_{\delta_2}} &\leq& c e^{-\epsilon |i-j|}
\|v_i\|_{L^2_{\delta_2}} \\ &\leq& c e^{-\epsilon |i-j|}
\|v_i\|_{L^p_{\delta}}
\end{eqnarray*}
now, note $\kappa_{ij}=1$ if $|i-j|\leq 1$ and $0$ otherwise;
using (\ref{est1}), we deduce
\begin{eqnarray*}
\|u_{ij}\|_{L^p_\delta} &\leq& c \big(
   \kappa_{ij}\|v_i\|_{L^p_\delta}
 + e^{-\epsilon |i-j|} \|v_i\|_{L^2_{\delta_2}} \big) \\
&\leq& c e^{-\epsilon |i-j|} \|v_i\|_{L^p_\delta}
\end{eqnarray*}
from which we deduce immediately $$ \|u\|_{L^p_\delta} \leq c
\|v\|_{L^p_\delta} , $$ which proves, with the help of local
elliptic regularity, that the operator is an isomorphism
$L^{2,p}_\delta \rightarrow L^p_\delta$. \pfend

\begin{rem}\label{neumann-cst}
Actually, the proof gives a bit more, namely the norm of the
inverse operator is bounded by a constant which is independent of
$R$ ($R$ big enough); this is because we have explicit constants
for the $L^2$ inverse, and the constants in the above proof do not
depend on $R$.
\end{rem}

\begin{rem}\label{neumann-holder}
The same proof works in H\"older spaces, and gives an isomorphism
between H\"older weighted spaces. In $C^k$ spaces, we have no more
elliptic regularity; nevertheless, if $v$ is in $C^0_\delta$, one
can still deduce from the above proof the estimate
\begin{equation}\label{estC0e}
\|r^{\delta-\epsilon} u\|_{C^0} \leq c \|r^\delta v\|_{C^0} ;
\end{equation}
this estimate is not a consequence of the $L^p$ estimate, because
the Sobolev embedding (which can be proven like the elliptic
estimate (\ref{est1}) by a homothety argument),
\begin{equation}\label{Sobolev-embedding}
\|u\|_{C^0_\delta} \leq c \big(
\|u\|_{L^p_{\delta-2/p}} + \|\nabla u\|_{L^p_{\delta+1-2/p}} \big),
\end{equation}
implies $L^{1,p}_{\delta+1-2/p} \hookrightarrow C^0_{\delta}$, so that
there is a loss of weight, since $L^{1,p}_{\delta+1-2/p}$ corresponds
to functions $O(r^{-\delta-1})$ when $C^0_{\delta}$ corresponds to functions
$O(r^{-\delta})$. Note also that in the case where $v$ lies in the
component where $\gamma$ acts non trivially, the maximum principle
provides directly the estimate (\ref{estC0e}) without the
$\epsilon$.
\end{rem}

\vskip12pt \centerline{\bf Analysis on 1-forms} \vskip10pt

In the next few lemmas, we assume that $a$ is a 1-form with values
in $End(E)$ such that $\partial_r\lrcorner a  = 0$ on $r=R$. Again
we suppose that $a$ is reduced to its component $a_\perp$. All
Sobolev norms are taken over the set $T\times V_R=\{r\geq R\}$.

\begin{lem} \label{l1}
One has the identity
\begin{equation}
\| d_{\Gamma}^*a \|_{L^2}^2 + \| d_{\Gamma}a \|_{L^2}^2 = \|
\nabla_{\Gamma}a \|_{L^2}^2 - \int_{r=R} \left| \frac{1}{r}
\frac{\partial}{\partial\theta}\lrcorner a \right|^2 dxdyd\theta
\end{equation}
with respect to the Euclidean metric.
\end{lem}

\pf The equality follows from the Weitzenb\"ock formula  in the
Euclidean metric: $$ d_{\Gamma}^*d_{\Gamma}+d_{\Gamma}d_{\Gamma}^*
= \nabla_\Gamma^*\nabla_\Gamma $$ Just integrate by parts and
check the boundary terms. \pfend

\begin{lem} \label{l3}
For any real function $f$ and any $R>0$, one has:
\begin{equation}
f(R)^2 \leq \frac{2}{R} \int_R^{R+1} (|\partial_r f|^2 + |f|^2) r
dr
\end{equation} \end{lem}

The proof is left to the reader.

\begin{lem} \label{l4}
If $R$ is sufficiently large, then for some constant $c$:
\begin{eqnarray*}
\| d_{\Gamma}^*a \|_{L^p_\delta} + \| d_{\Gamma}a \|_{L^p_\delta}
&\geq& c \| \nabla_{\Gamma}a \|_{L^p_\delta} \\ \| d_{\Gamma}^*a
\|_{C^0_\delta} + \| d_{\Gamma}a \|_{C^0_\delta} &\geq& c \| a
\|_{C^0_{\delta-\epsilon}}
\end{eqnarray*}
with respect to the Euclidean metric.
\end{lem}

\begin{rem} Remind that on the component $a=a_\perp$ we look at,
$\nabla_\Gamma a$ controls $a$ by (\ref{est0}).
\end{rem}

\begin{proof}
From lemma \ref{l3} and lemma \ref{l2}, we have:
\begin{eqnarray*}
\int_{r=R} | a |^2 dx dy & \leq & \frac{2}{R} \int_{[R,R+1]}
                                (|\nabla_{\partial_r} a|^2 +
                                |a|^2) r dr dx dy \\
&\leq & \frac{C_1}{R} \int_{[R,R+1]} |\nabla_\Gamma a|^2 dx dy r
dr
\end{eqnarray*}
for some constant $C_1$; in particular $$ \int_{r=R} \left|
\frac{1}{r} \frac{\partial}{\partial\theta}\lrcorner a \right|^2
dxdyd\theta \leq \frac{C_1}{R} \int_{[R,R+1]} |\nabla_\Gamma a|^2
dx dy r dr d\theta $$ and we deduce from lemma \ref{l1}, for $R$
big enough, $$ \| d_{\Gamma}^*a \|_{L^2}^2 + \| d_{\Gamma}a
\|_{L^2}^2 \geq \frac{1}{2} \| \nabla_{\Gamma}a \|_{L^2}^2 $$
which proves the $L^2$-estimate of the lemma.

The $L^2$-estimate with weights is proven in the same way. In the
integration by parts, new terms appear because of the weight
$r^{2\delta}$. However, as in the proof of lemma \ref{neumann-L2},
these terms have all a coefficient $O(r^{-1})$ and therefore are a
small perturbation if $R$ is large enough (note that we can take
the same $R$ if the weight remains bounded).

Finally, one may deduce the $L^p$ and $C^0$ estimates from the
$L^2$ estimates as in lemma \ref{neumann-Lp} and remark
\ref{neumann-holder}, since the operator $d_\Gamma^*+d_\Gamma$ has
injective symbol, and the boundary condition $\partial_r\lrcorner
a=0$ is an elliptic boundary condition. The proof is a slightly
more complicated, because one has to compose the decomposition
(\ref{vi}) with a $L^2$-projection on the image of the operator.
\end{proof}

\begin{lem} \label{defn.op}
The operator $2 d_\Gamma^*d_\Gamma^+ + d_\Gamma d_\Gamma^*$ on
1-forms lying in $\Omega^1\otimes End(E)_\perp$, with Dirichlet
condition on $r=R$, is an isomorphism in weighted Sobolev or
H\"older spaces for all weights $\delta$.
\end{lem}

\begin{proof} Again the Weitzenb\"ock formula
$$ 2 d^*_\Gamma d_\Gamma^+ + d_\Gamma d_\Gamma^* = \nabla_\Gamma^*
\nabla_\Gamma $$ gives the $L^2$-estimate (for forms vanishing on
the boundary)
\begin{eqnarray*}
((2 d^*_\Gamma d_\Gamma^+ + d_\Gamma d_\Gamma^*) u,u)_{L^2} &=&
\|\nabla_\Gamma u\|_{L^2}^2 \\ &\geq& c \|u\|_{L^2}^2
\end{eqnarray*}
from which the $L^2$-statement (without weight) follows
immediately. One can then deduce weighted statements as in the
proofs of lemmas \ref{neumann-L2} and \ref{neumann-Lp}.
\end{proof}

%-----------------------------------------------------------------------------

\section{Existence of a Coulomb gauge} \label{coulomb}

After the technical work of the previous section, we are finally
in a position to establish theorem \ref{coulomb.gauge}, the key
analytical result of this paper. The first step is the nonlinear
version of the H\"older estimate in lemma \ref{l4}; the exponent
$p$ is fixed, near infinity.

\begin{lem} \label{l9}
Given $\eta_1$ sufficiently small, if a connection
$A=\Gamma+a_\Gamma+a_\perp$ on $r\geq R$ satisfies:
\begin{enumerate}
\item $d^*_{\Gamma+a_\Gamma} a_\perp = 0$ ,
\item $\partial_r \lrcorner a (r=R) = 0$ ,
\item $\| r^\epsilon a \|_{C^0} \leq \eta_1$ ,
\end{enumerate}
then:
\begin{equation} \label{l9eqn}
\| r^2 F_{\Gamma+a_\Gamma} \|_{C^0} + \| r^{2-\epsilon} a_\perp
\|_{C^0} +
\|(\nabla_\Gamma+a_\Gamma)a_\perp\|_{L^p_{2-2/p-\epsilon}} \leq c
\|r^2 F_A \|_{C^0} .
\end{equation} \end{lem}

\begin{proof} First, note that:
\begin{equation}
F_A = F_{\Gamma+a_\Gamma} + d_{\Gamma+a_\Gamma} a_\perp +
\frac{1}{2}[a_\perp,a_\perp] .
\end{equation}
Therefore, using the decomposition in (\ref{decomposition01}), we
have:
\begin{eqnarray}
\big( F_A \big)_\Gamma & = & F_{\Gamma+a_\Gamma}
  + \frac{1}{2} \big( [a_\perp,a_\perp] \big)_\Gamma \label{decomp.curv1} , \\
\big( F_A \big)_\perp & = & d_{\Gamma+a_\Gamma}a_\perp
  + \frac{1}{2} \big( [a_\perp,a_\perp] \big)_\perp , \label{decomp.curv2}
\end{eqnarray}
from which the the estimates below follow:
\begin{eqnarray}
\| r^2 \big( F_A \big)_\Gamma \|_{C^0} & \geq & \| r^2
F_{\Gamma+a_\Gamma} \|_{C^0}
  - \| r a_\perp \|_{C^0}^2 \label{decomp.est1} , \\
\| r^2 \big( F_A \big)_\perp \|_{C^0} & \geq & \| r^2 d_\Gamma
a_\perp \|_{C^0}
  - \| r^2 [a_\Gamma,a_\perp] \|_{C^0} - \| r a_\perp \|_{C^0}^2 .\label{decomp.est2}
\end{eqnarray}
Using $C^0_2\subset L^p_{2-2/p-\epsilon}$ and the estimate in
lemma \ref{l4}, we get:
\begin{eqnarray*}
\|r^2 F_A \|_{C^0} & \geq & c \bigg( \| r^2 F_{\Gamma+a_\Gamma}
\|_{C^0} + \| r^{2-\epsilon} a_\perp \|_{C^0} +
\|(\nabla_\Gamma+a_\Gamma)a_\perp\|_{L^p_{2-2/p-\epsilon}} \bigg)
\\ & & \quad - c' \bigg(\| r a_\perp \|_{C^0}^2 + \| r^2
[a_\Gamma,a_\perp] \|_{C^0}\bigg);
\end{eqnarray*}
from the third hypothesis, we have $$ \| r a_\perp \|_{C^0}^2 + \|
r^2 [a_\Gamma,a_\perp] \|_{C^0} \leq \eta_1 \| r^{2-\epsilon}
a_\perp \|_{C^0} ; $$ if $\eta_1$ is small enough, these two
inequalities give the required estimate.
\end{proof}

\begin{lem} \label{cpt.sup}
Given $\eta$, there exists $R$ such that if $A$ is a connection
over $T\times V_R$ such that $A-\Gamma$ is compactly supported and
$|F_A| \leq \eta\cdot r^{-2}$, then there is a gauge $g$ such that
$g(A)=\nabla_\Gamma+a_\Gamma+a_\perp$, with:
\begin{eqnarray*}
\textnormal{(i)} & & d^*_{\Gamma+a_\Gamma} a_\perp = 0, \\
\textnormal{(ii)} & & \partial_r \lrcorner a (r=R) = 0, \\
\textnormal{(iii)} & & \| r^2 F_{\Gamma+a_\Gamma} \|_{C^0}
  + \| r^{2-\epsilon} a_\perp \|_{C^0}
  + \|(\nabla_\Gamma+a_\Gamma)a_\perp\|_{L^p_{2-2/p-\epsilon}}
  \leq c \| r^2 F_A \|_{C^0} .
\end{eqnarray*}\end{lem}
\pf We now have all the necessary ingredients for a proof by
continuity. Consider the homothety $\phi_t(r)=e^t r$ and the
connections $A_t=\phi_t^*A$. We have $A_0=A$ and, for $t$ big
enough, say $t\geq T$, $A_t=d_\Gamma$ because of the assumption on
compact support. Moreover, it is clear from the form of the metric
that $$ |F_{A_t}| = |\phi_t^* F_A| \leq \phi_t^* |F_A| \leq
\frac{c e^{-2t}}{r^2} $$ so that the whole path of connections
$(A_t)$ satisfies the hypothesis of the lemma. Moreover, after
gauge transformation, we can also assume that $A_t=\Gamma+a_t$
with $\partial_r \lrcorner a_t (r=R) = 0$ for all $t$.

We prove that the subset $S\subseteq[0,T]$ containing all the
values of $t$ for which the theorem holds for $A_t$ is both closed
and open. Since $S$ is nonempty (it contains $t=T$), $S$ must be
the whole interval and the result holds for $t=0$.

The closedness is trivial, since the estimate on the connection
provides all the needed bounds.

For openness, first remark that proposition
\ref{invariant-control} provides a gauge in which
\begin{eqnarray*}
\| r^2 F_{d_\Gamma+a_\Gamma} \|_{C^0} &\geq& c \| \frac{r}{\ln r}
a_\Gamma \|_{C^0} \\ &\geq& c \frac{R^{1-\epsilon}}{\ln R} \|
r^\epsilon a_\Gamma \|_{C^0} ;
\end{eqnarray*}
on the other hand, from (iii),
\begin{eqnarray*}
\| r^2 F_A \|_{C^0} &\geq& c \|r^{2-\epsilon} a_\perp\|_{C^0} \\
&\geq& c R^{2-2\epsilon} \|r^\epsilon a_\perp\|_{C^0} ;
\end{eqnarray*}
we deduce
\begin{equation}
\|r^\epsilon a\|_{C^0} \leq c^{-1} R^{-(1-2\epsilon)} \eta ;
\end{equation}
taking $R$ big enough so that the RHS is smaller than $\eta_1$ of
lemma \ref{l9}, we see that (i) and (ii) imply (iii).

It remains to solve problem (i)-(ii) near a solution. Fix some $t$
and suppose that $g_t(A_t)=\Gamma+b$ with $\Gamma+b$ satisfying
(i), (ii) and (iii). If we have a connection  $\Gamma+b+\varpi$
with $\partial_r \lrcorner \varpi (r=R) = 0$, we want to find a
gauge $g$ such that: $$ \left\{ \begin{array}{l} g(\Gamma + b +
\varpi) = \Gamma + c_\Gamma + c_\perp \\ d^*_{\Gamma+c_\Gamma}
c_\perp = 0
\end{array} \right. $$
Looking at solutions of the form $g=e^{u_\perp}$, the equation to
be solved is: $$ L(u_\perp,\varpi)= d^*_{\Gamma+c_\Gamma} \big(
e^u(\Gamma+c_\Gamma+c_\perp) e^{-u} -
d_{\Gamma+c_\Gamma}(e^u)\cdot e^{-u} \big) = 0 ; $$ we would like
to solve this equation with $u_\perp$ in $C^2$, but $C^k$ spaces
are not suitable for elliptic analysis; instead, we use weighted
$L^p$ spaces with $p$ very big; since we have the freedom to apply
a $\Gamma$-invariant gauge transformation, using proposition
\ref{invariant-control}, we can choose a gauge in which the
derivatives of $b_\Gamma$ are also controled, and therefore the
operator $L$ is well defined; its linearization along the first
variable is given by the operator: $$ u \ \seta \ d_\Gamma^*
d_\Gamma u + \textrm{perturbation}; $$ if $R$ is big enough, the
perturbation is sufficiently small and we get an isomorphism by
lemma \ref{neumann-L2}.

This completes the proof. \pfend

\paragraph{Completing the proof of theorem \ref{coulomb.gauge}.}
Our final task is to remove from lemma \ref{cpt.sup} the
assumption that $A-\Gamma$ is compactly supported.

Using proposition \ref{flat-limit}, we approximate the connection
$A$ by a sequence $A_i$ such that $\Gamma-A_i$ is compactly
supported, and $\|r^2 F_{A_i}\|_{C^0}$ remains bounded.

We can apply lemma \ref{cpt.sup} to each connection $A_i$, thus
obtaining a gauge $g_i$ such that $g_i(A_i)=d_\Gamma+a_i$, and
$a_i$ satisfies (i)--(iii) of lemma \ref{cpt.sup}. Using
proposition \ref{invariant-control} for the $(a_i)_\Gamma$ part,
the $(a_i)$ converge (weakly) to a limit $a$ still satisfying
(i)--(iii), such that $d_\Gamma+a$ is gauge equivalent to $A$.
\pfend

\part{Instantons, holomorphic bundles, and the moduli space} \label{part2}

So far, $A$ has simply been a connection on $E\seta\torus$ with 
quadratic curvature decay. From now on, we shall assume that $A$ 
is also an instanton.

\section{Asymptotic behavior: proof of theorem \ref{asymp.par}} \label{ab}

Let us now assume that $A$ is a doubly-periodic instanton
connection. Using theorem \ref{coulomb.gauge}, if $R$ is big
enough, we can put it in a Coulomb gauge on $r\geq R$, so that
$A=\Gamma+a_\Gamma+a_\perp$, with $a_\Gamma$ and $a_\perp$
satisfying the Coulomb gauge equation, $$ d_{\Gamma+a_\Gamma}^*
a_\perp =  0 , $$ and the anti-self-duality equation, $$ d^+_\Gamma
a + \frac{1}{2}[a,a]^+ = 0 . $$ These can be rewritten as follows:
\begin{eqnarray}
d_\Gamma^* a_\perp &=& -a_\Gamma^*a_\perp \label{eq1} \\
 d_\Gamma^+ a_\perp &=& -[a_\Gamma,a_\perp]^+ -
\frac{1}{2}[a_\perp,a_\perp]^+_\perp \label{eq2} \\
 d^+a_\Gamma + \frac{1}{2} [a_\Gamma,a_\Gamma]^+ &=&
  -[a_\perp,a_\perp]^+_\Gamma \label{eq3}
\end{eqnarray}

Now let $\chi=\chi(r)$ be a smooth cut-off function supported on
$T\times V_R$; we have, using equations (\ref{eq2}) and
(\ref{eq3}):
\begin{equation} \label{defn.id}
\big( d_\Gamma^+ + d_\Gamma^* \big) (\chi a_\perp) = \chi
(a_\Gamma\odot a_\perp + a_\perp\odot a_\perp) + d\chi\odot
a_\perp
\end{equation}
where $\odot$ denotes some bilinear operations.

From theorem \ref{coulomb.gauge} and proposition
\ref{invariant-control}, we already know that
$|a_\perp|=O(r^{-2+\epsilon})$ and that we can choose a gauge such
that $|a_\Gamma|=O(\ln r/r)$. We now apply lemma \ref{defn.op} to
the equation (\ref{defn.id}): a priori the lemma applies to the
laplacian $(d_\Gamma^+)^* d_\Gamma^+ + d_\Gamma d_\Gamma^*$ but
the estimates also imply estimates for the first order elliptic
operator $d_\Gamma^+ + d_\Gamma^*$ (alternatively one may take one
derivative of equation (\ref{defn.id}) and use the bounds on the
derivatives of $a_\perp$ and $a_\Gamma$); the RHS of equation
(\ref{defn.id}) is $O(r^{-3+\epsilon})$, therefore $|a_\perp|=
O(r^{-3+\epsilon_2})$, where $\epsilon_2>\epsilon$; by the same
argument, we have that $|a_\perp|=O(r^{-4+\epsilon_3})$, etc.
Therefore, $|a_\perp|=O(r^{-\delta})$ for any $\delta>0$.

Now come back to equation (\ref{eq1}): it now means that
$d+a_\Gamma$ satisfies the instanton equation up to a term which
goes very quickly to $0$ at infinity; as $a_\Gamma$ is translation
invariant, this means, by dimensional reduction, that $d+a_\Gamma$
is a solution of Hitchin's equations for Higgs bundles on $\bR^2$
near infinity, up to a term decaying quicker than any
$O(r^{-\delta})$. The behavior of the solutions of Hitchin's
equations near a singularity has been studied by Simpson \cite{S},
Biquard \cite{B3}. The arguments in these papers are not affected
by a very quickly decaying perturbation. Moreover, the bounds in
proposition \ref{invariant-control} implies that the Higgs field
is $O(1/r)$ at infinity, so that the Higgs bundle is ``tame'' in
Simpson's terminology. Finally, we deduce from these articles that
$d+a_\Gamma$ is close to one of the examples described in the
introduction, in the sense of theorem \ref{asymp.par}. \qed

%-----------------------------------------------------------------------------

\section{Holomorphic extension} \label{holo.extn}

The theorem \ref{asymp.par} proves that any instanton $A$ with
quadratic curvature decay can be put in a gauge near infinity so that
$$ A = A_0 + a , $$
where $A_0$ is one of the model torus invariant instantons
induced by model Higgs bundles, and $a$ is a small perturbation.

\vskip12pt \centerline{\textbf{Local aspects}} \vskip10pt

Let us now restrict to the semisimple case. Therefore, we have
$$ A_0=d+i\left(\begin{array}{cc}a_0 & 0 \\ 0 & - a_0 \end{array}\right) $$
with
$$ a_0 = \lambda_1 dx + \lambda_2 dy + (\mu_1 \cos\theta - \mu_2 \sin\theta)
\frac{dx}{r} + (\mu_1 \sin\theta + \mu_2\cos\theta) \frac{dy}{r} + \alpha d\theta ; $$
observe that the (0,1)-part of this form is
$$
a_0^{0,1} = \lambda d\overline{z} +
            \mu \frac{d\overline{z}}{w} -
            \frac{\alpha}{2} \frac{d\overline{w}}{\overline{w}} ,
            \quad
\lambda=\frac{\lambda_1+i\lambda_2}{2}, \, \mu=\frac{\mu_1+i\mu_2}{2},
$$
so there is a singularity in the direction of transverse disks to the
torus at infinity. We first reduce to a normal form on transverse disks.

\begin{lem}
Near the torus at infinity, there exists a continuous complex gauge
transformation $g$, such that
\begin{enumerate}
\item $g|_{T_\infty}=1$ ;
\item $|\nabla_{A_0}g g^{-1}|=O(r^{-(1+\delta)})$ (and $g_\perp$ is
$O(r^{-\delta})$ for any $\delta$);
\item $g(\overline{\partial}_{A}) = A_0 + b d\overline{z}$, with
$b=O(r^{-(1+\delta)})$.
\end{enumerate}
\end{lem}

\begin{proof}
We give a concise proof, since this is parallel to \cite[section 9]{B3}.
Remark that
\begin{equation}\label{dbalpha}
\overline{\partial}_\alpha = \overline{\partial}
  - \frac{\alpha}{2} \frac{d\overline{w}}{\overline{w}} =
  r^{-\alpha} \circ \overline{\partial} \circ r^\alpha ;
\end{equation}
now the problem to be solved is
$$\frac{\overline{\partial}_{A_0}}{\partial\overline{w}} g - ga = 0 , $$
that is, using $g=1+u$,
$$ \bigg( \frac{\partial}{\partial\overline{w}} - \frac{1}{2} \left(
\begin{array}{cc} \alpha & 0 \\ 0 & -\alpha \end{array} \right) \bigg)
u - u a = -a ; $$
this is a $\overline{\partial}$-problem on small disks near infinity; for the
model problem (\ref{dbalpha}) the Cauchy formula gives us an explicit
solution; in general, with the small perturbation $a$, the solution is
produced by a fixed point theorem, and we even have an estimate
$$ \sup r^\delta |u| \leq c \sup r^{1+\delta}|a| ; $$
one can then deduce the regularity statement on $u$.
\end{proof}

Note $b_{jk}$ the coefficients of the matrix $b$ above.
Let $(e_1,e_2)$ be the orthonormal basis for the trivialisation of the
bundle near infinity. From the lemma and equation (\ref{dbalpha}), we
deduce that the sections
\begin{equation}\label{basis-hol-ext}
(\sigma_1=r^{-\alpha}g(e_1),\sigma_2=r^\alpha g(e_2))
\end{equation}
are holomorphic on transverse disks, and, moreover, in the basis
$(\sigma_1,\sigma_2)$, we now have
\begin{equation}\label{dbsigma}
\overline{\partial}_A = \overline{\partial}
 + \left(\begin{array}{cc} \lambda & 0 \\ 0 & -\lambda
         \end{array}\right) d\overline{z}
 + \left(\begin{array}{cc} \mu & 0 \\ 0 & -\mu
         \end{array}\right) \frac{d\overline{z}}{w}
 + \left(\begin{array}{cc} b_{11} & r^{2\alpha}b_{12} \\
         r^{-2\alpha}b_{21} & b_{22} \end{array}\right) d\overline{z} , 
\end{equation}
with all coefficients of the last matrix holomorphic in $w$.
From this, we see immediately that in the basis $(\sigma_1,\sigma_2)$, 
the operator (\ref{dbsigma}) defines a holomorphic extension $\ee$ over 
$\tproj$.

Since
$$ |\sigma_1| \sim r^{-\alpha} , \quad |\sigma_2| \sim r^{\alpha} , $$
we see that, from an intrinsic point of view, if $\alpha<1/2$, the local holomorphic sections
of $\ee$ are characterized as the local holomorphic sections $\sigma$
outside $T_\infty$ satisfying the growth condition
\begin{equation}\label{pol-growth}
|\sigma| = O(r^{\alpha}) .
\end{equation}

When $0<\alpha<1/2$, this global
extension has a subbundle $\ff$ over the torus at infinity, given by the
values of the local holomorphic sections $\sigma$ satisfying the
growth condition
\begin{equation}\label{pol-sub-growth}
|\sigma| = O(r^{-\alpha}) .
\end{equation}
Therefore, the growth of the holomorphic sections at infinity
determine a ``parabolic structure''
$$ \ee \supset \ff \supset 0 , $$
with weights $-\alpha < \alpha$ (the sign is changed because the local
coordinate near infinity is $w^{-1}$).

Actually one can say more : over $T_\infty$, the 
$\overline{\partial}$-operator (\ref{dbsigma}) is
$$ \overline{\partial} +
  \left(\begin{array}{cc} \lambda & 0 \\ 0 & -\lambda
	\end{array} \right) d\overline{z} , $$
which means that
$$ \ee|_{T_\infty} = L_{\xi_0} \oplus L_{-\xi_0} . $$
Of course, if $\alpha$ is nontrivial,
then $\ff=L_{\xi_0}$ is canonically determined by the growth
condition (\ref{pol-sub-growth}).

Actually, the decomposition $L_{\xi_0}\oplus L_{-\xi_0}$ can almost always be
made canonical: this is clear if $\xi_0\neq 0$, and in this case,
since the off-diagonal components of the connection decay quicker than any
$O(r^{-\delta})$, we deduce from equation (\ref{dbsigma}) that, still
in the basis $(\sigma_1,\sigma_2)$,
\begin{equation}\label{dev-ee}
\overline{\partial}_A = \overline{\partial}
 + \left(\begin{array}{cc} \lambda & 0 \\ 0 & -\lambda
         \end{array}\right) d\overline{z}
 + \left(\begin{array}{cc} \mu & 0 \\ 0 & -\mu
         \end{array}\right) \frac{d\overline{z}}{w}
    + O(r^{-2}) ;
\end{equation}
this gives the asymptotic behavior of $\ee|_{T_w}$ when $w$ goes to infinity.

Moreover, when $\xi_0=0$, we still get something from (\ref{dbsigma}):
since the coefficients are holomorphic in $w$, we note $b'_{12}$ the
coefficient of $r^{2\alpha}b_{12}$ on $w^{-1}$ (in the case $\alpha=0$, we
simply have $b'_{12}=0$), so
\begin{equation}\label{dev-ee2}
\overline{\partial}_A = \overline{\partial}
 + \left(\begin{array}{cc} \mu & b'_{12} \\ 0 & -\mu
         \end{array}\right) \frac{d\overline{z}}{w}
 + O(r^{-2}) ;
\end{equation}
if $\mu\neq 0$, the matrix appearing above can always been
diagonalized with eigenvalues $\pm\mu$, which means that up to
changing $\sigma_2$ by some multiple of $\sigma_1$, we are reduced to
(\ref{dev-ee}) so that a supplementary subspace of $\ff$ is still well
defined (and when $\alpha=0$, the decomposition $\cpx\oplus\cpx$ still
makes sense, as the eigenspaces of this matrix).

Note also that, as a consequence of (\ref{basis-hol-ext}), since $g$
is continuous, the unitary extension (given by the basis $(e_1,e_2)$
of the Coulomb gauge) and the holomorphic extension are topologically
isomorphic.

Therefore, we have proven the following proposition.

\begin{prop}\label{prop-extn}
In the semisimple case, for $\alpha<1/2$,
if $A$ is a doubly-periodic instanton connection satisfying
$|F_A| = O(r^{-2})$, then $A^{0,1}$ has a unique holomorphic
extension $\ee$ over $\tproj$, whose holomorphic sections satisfy the growth
condition (\ref{pol-growth}). Moreover, one has $c_2(\ee)=k$ and a
decomposition (if $\lambda$ or $\mu$ is nonzero)
$\ee|_{T_\infty} = L_{\xi_0} \oplus L_{-\xi_0}$.\qed
\end{prop}

\begin{rem}
Note that when $\alpha=1/2$, we cannot get a $Sl_2$-extension this way:
indeed we could equally well choose the sections
$(w\sigma_1,\sigma_2/w)$, giving a different extension.
One way to construct a canonical extension is to use
(\ref{pol-growth}) with $\alpha=-1/2$, which furnishes a
$Gl_2$-extension where all nonzero sections have norm $O(r^{-1/2})$.
Also, a $Sl_2$-extension can be constructed if $\xi_0\neq-\xi_0$,
by \emph{deciding} that sections with nonzero values in
$L_{\pm\xi_0}$ have norm like $r^{\mp 1/2}$.

In the sequel we will ignore this case, but all the statements can be
easily adapted to it.
\end{rem}

\begin{rem}
In the nilpotent case (then $\lambda$, $\mu$, and $\alpha$ are  trivial),
the result is the same, but (as in the case of Higgs bundles) the
growth of the holomorphic sections at infinity is now logarithmic:
\begin{equation}\label{ln-growth}
|\sigma| = O\big((\ln r)^{\frac{1}{2}}\big) ,
\end{equation}
and there is a line subbundle $\ff$ defined by the growth condition
\begin{equation}\label{ln-sub-growth}
|\sigma| = O\big((\ln r)^{-\frac{1}{2}}\big) . 
\end{equation}
The subbundle $\ff$ has no canonical supplementary subspace.
The tools in \cite[section 9]{B3} handle this situation as well.

Also observe that the $\del$-operator for the model instanton
(\ref{ex3}) is (in an orthonormal basis $(e_1,e_2)$)
$$ \del + \left( \begin{array}{cc}1&0\\0&-1\end{array}
          \right) \frac{d\overline{w}}{2\overline{w}\ln r^2}
        + \frac{1}{r\ln r^2}
  \left( \begin{array}{cc}0&e^{-i\theta}d\overline{z}\\0&0\end{array}
  \right)
$$
which gives, in the basis
$(e_1/(\ln r^2)^{\frac{1}{2}},e_2 (\ln r^2)^{\frac{1}{2}})$,
$$
\del + \left( \begin{array}{cc}0&\frac{d\overline{z}}{w}\\0&0\end{array}
       \right) ;
$$
in particular, $\ee|_{T_w}$ is the nontrivial extension of
$\underline{\cpx}$ by $\underline{\cpx}$; it is easy to see that this
remains true for instantons, asymptotic to this nilpotent model.
\end{rem}

\vskip12pt \centerline{\textbf{Non-existence results}} \vskip10pt

The proposition \ref{prop-extn} gives obstructions for the existence
of instantons. Here are some examples.

\begin{lem} \label{non-ex1}
There are no instantons with $\xi_0=-\xi_0$ and $k=1$.
\end{lem}
\begin{proof}
For a contradiction, let $A$ be an instanton with $\xi_0=-\xi_0$ and $k=1,2$,
and consider the extended holomorphic bundle $\ee$ given by theorem \ref{extn}.
The restriction of $\ee$ to the elliptic fibres $T_p$ must be semistable for
all $p\in\proj$ (see \cite{J3}). Moreover, $\ee|_{T_p}$ cannot be generically 
the nontrivial extension of $\underline{\cpx}$ by itself, since this would give
a non-constant map from $\proj$ to $\cpx$ (which parametrises the extensions of 
$\underline{\cpx}$ by itself).

Therefore, as shown in \cite{J2,J3}, index theory tells us that for each $\xi\in\dual$:
\begin{equation} 
\Sigma_{w\in\proj} h^0 \big( T_w,\ee\otimes L_\xi|_{T_w} \big) = k
\end{equation}
But if $\ee|_{T_\infty}=L_{\xi_0}\oplus L_{\xi_0}$, then
$ h^0 \big( T_\infty,\ee\otimes L_{\xi_0}|_{T_\infty} \big) = 2 $, thus 
contradicting the assumption that $k=1$.
\end{proof}

\begin{lem} \label{non-ex2}
There are no instantons with $\xi_0\neq-\xi_0$ and $\mu=0$.
\end{lem}
\begin{proof}
The lemma is a consequence of the Nahm transform of doubly-periodic instantons
defined in \cite{J2}, more exactly of its holomorphic aspects; we
anticipate a bit here, but see the introduction to Part III for a summary
of the construction.

Again for a contradiction, let $A$ be an instanton with $\mu=0$ and 
asymptotic state $\as$ not of order two. The corresponding Nahm 
transformed Higgs field $\Phi$ has simple poles at $\as$; its residues 
have rank one. However, as we shall see in the proof of theorem 
\ref{ap.nt.thm}, the non-zero eigenvalues of the residues of $\Phi$ 
are exactly $\pm\mu$, and more generally, the eigenvalues of $\Phi$ at 
$\xi\in\dual$ are the $w$ such that $H^0(T_w,\ee\otimes L_{\xi})\neq 0$; 
hence, the vanishing of $\mu$ implies that the eigenvalues of $\Phi$ remain
bounded when $\xi$ goes to $\xi_0$.

Now if $\xi_0\neq\xi_0$ then $\ee$ remains isomorphic to some 
$L_\xi\oplus L_{-\xi}$ on each torus near infinity. It is then 
clear (again, see the proof of theorem \ref{ap.nt.thm}) that 
the eigenvalues of $\Phi$ must go to infinity and we get a
contradiction.
\end{proof}

\vskip12pt \centerline{\textbf{Global aspects, stability}} \vskip10pt

More subtil obstructions come from stability properties.
We investigate this for the extension $\ee$ of an instanton $A$ with
quadratic curvature decay. Notice that by theorem
\ref{asymp.par}, in the semisimple case, the curvature is only $O(r^{-2})$, but
\begin{equation}\label{cont-curv-ss}
|\iota_{\{\cdot\}\times \cpx}F_A| + |\iota_{T \times \{\cdot\}}F_A|
 = O\big(r^{-(2+\epsilon)}\big) ;
\end{equation}
in the nilpotent case, we have
\begin{equation}\label{cont-curv-n}
|F_A| = O\big(r^{-2} (\ln r)^{-2}\big) ;
\end{equation}
the point here is that these two controlling factors are in $L^1$, whence $F_A$
itself is not $L^1$: this will enable us to define a degree.

The degree of a saturated subsheaf $L$ of $\ee$ with respect to the Euclidean
K\"ahler form $\omega$ is \cite[lemma 3.2]{S0}
\begin{equation}\label{for-deg}
2\pi \deg L = i \int \textnormal{tr}(\pi F_{A})\wedge \omega
            - \int |\overline{\partial}\pi|^2
\end{equation}
where $\pi$ is orthogonal projection on $L$; from (\ref{cont-curv-ss})
and (\ref{cont-curv-n}), this can be $-\infty$ or a real number; in
the last case, $\overline{\partial}\pi$ is in $L^2$: this condition
must be analyzed more precisely.

Again, we now restrict to the semisimple case (see remark
\ref{sta-nil} for the nilpotent case),
so that $\ee|_{T_\infty}=L_{\xi_0}\oplus L_{-\xi_0}$,
with weights $-\alpha$ and $\alpha$, and behavior (\ref{dev-ee}).
In this case, we have near infinity
\begin{equation}\label{sum-infin}
\ee|_{T_w} = L_{\xi(w)} \oplus L_{-\xi(w)} .
\end{equation}
\begin{lem}\label{finite-deg}
Suppose $\alpha\neq 0$, then
the degree of a subsheaf $L$ of $\ee$ is finite if and only if
\begin{enumerate}
\item $L|_{T_\infty}$ is flat; in particular, if $L_{\xi_0}\neq L_{-\xi_0}$,
this means that $L\subset L_{\pm\xi_0}$ ;
\item if $L|_{T_\infty} \subset \ff=L_{\xi_0}$, then
$L|_{T_w} \subset L_{\xi(w)}$ up to first order near infinity.
\end{enumerate}
Now suppose $\alpha=0$, then the degree of a subsheaf $L$ of $\ee$ is
finite if and only if $L|_{T_w} \subset L_{\pm\xi(w)}$ up to first
order near infinity.
\end{lem}

\begin{rem}
The first order condition can be seen as a reminiscence of the
approximating Higgs bundle at infinity; indeed the Higgs field has
eigenspaces $L_{\pm\xi(w)}$ and for Higgs bundle stability, one looks
only at subsheaves stable under the action of the Higgs field.
\end{rem}

\begin{proof}
We analyze the situation locally near infinity; in the decomposition
(\ref{sum-infin}), the metric is approximately
$$ \left(\begin{array}{cc} r^{-2\alpha} & 0 \\ 0 & r^{2\alpha}
	 \end{array}\right), $$
and we will simplify the problem by using this metric to make the
calculations (the correction term can be easily bounded);
at a point on $T_\infty$ where $L$ is a subbundle, we suppose for
example that $L$ is not contained in $L_{-\xi_0}$; choose a local flat
section $\sigma$ for $L_{\xi_0}$, and note $\sigma^t$ the dual flat
section of $L_{-\xi_0}$; extend $\sigma$ near $T_\infty$, keeping it
parallel on $T_w$ (this is possible with our approximation for the
metric); locally, $L$ is generated by $s=\sigma+f\sigma^t$, where $f$ is
holomorphic, and an orthogonal section is given by
$t=r^{2\alpha}\sigma-\overline{f}\sigma^t r^{-2\alpha}$,
and 
$$ \overline{\partial}_T t = - (\overline{\partial}_T \overline{f})
  \sigma^t r^{-2\alpha},$$
from which we deduce
$$ \pi(\overline{\partial}_Tt) = -
 \frac{\overline{\partial}_T\overline{f}}{r^{-2\alpha}+|f|^2r^{2\alpha}} s,
$$
and finally, since our choice of $t$ satisfies $|s|=|t|$, and $f$ is
holomorphic, 
$$
|\overline{\partial}_T\pi| =
\frac{|d_T f|}{r^{-2\alpha}+|f|^2r^{2\alpha}};
$$
in order for $\overline{\partial}_T\pi$ to be in $L^2$, it is
necessary that $d_T f=0$ on $T_\infty$, and therefore $L$ is
constant.

Now restrict to the case of nontrivial decomposition $L_{\xi_0}\oplus
L_{-\xi_0}$ (the other cases are similar);
therefore we may suppose that $f=0$ on $T_\infty$; if the
first order term of $d_T f$ does not vanish, then
$$ |\overline{\partial}_T\pi| \sim r^{-1+2\alpha} $$
this still is not in $L^2$ if $\alpha\geq 0$ (but it is in $L^2$ if
$\alpha<0$, which corresponds to the case $L|_{T_\infty}\subset L_{-\xi_0}$);
this means that we need $d_T f$ to vanish up to first order.

Concerning $\overline{\partial}_{\cpx}\pi$, it is easy to verify that
the $L^2$-condition is always satisfied.
\end{proof}

Recall that
\begin{equation}\label{for-FL}
F_L = \pi F_A \pi + \overline{\partial}\pi \wedge \partial \pi .
\end{equation}
When the degree is finite, that is when $\overline{\partial}\pi$ is
$L^2$, the restriction of $\omega$ to $\cpx$ does not contribute:
indeed, $dw\wedge d\overline{w}=\partial\overline{\partial}|w|^2$,
and this leads to
$$
\int_{r\leq R} F_L \wedge dw\wedge d\overline{w} =
\int_{r=R} w dw \wedge F_L
$$
but using (\ref{cont-curv-ss}) and (\ref{for-FL}), we see that this
goes to zero as $R$ goes to infinity. Then we can rewrite the degree
(denoting $\overline{\partial}_{\cpx}$ the
$\overline{\partial}$ operator in the $\cpx$ direction)
\begin{equation}\label{forbis-deg}
2\pi \deg L =
i \int \pi F_{A}\wedge \omega_T - \int |\overline{\partial}_{\cpx}\pi|^2 ,
\end{equation}
and this in turn is easily interpreted \cite[(4.1)]{B2} as a
``parabolic degree'': 
\begin{equation}\label{par-deg}
\deg L = \left\{ \begin{array}{l}
  c_1(L)[t] + \alpha \langle[\omega_T],[t]\rangle \quad \mbox{ if }
    L_{T_\infty}\subset L_{-\xi_0}, \\
  c_1(L)[t] - \alpha \langle[\omega_T],[t]\rangle \quad \mbox{ if }
    L_{T_\infty}\subset L_{\xi_0}, \end{array} \right.
\end{equation}
where $[t]$ is the fundamental class of $T$ and $\omega_T$ the given
K\"ahler form on $T$; of course this is not a degree in the usual
sense on $\tproj$, since we use the non ample class $[t]$.

Define $\alpha$-\emph{stability} of $\ee$ as the fact that any
subsheaf satisfying the condition of lemma \ref{finite-deg} has
negative degree (we shall forget the $\alpha$ when there is no ambiguity); 
standard arguments give us

\begin{prop}\label{stable}
If $A$ is an instanton with quadratic curvature decay, then the
holomorphic extension $\ee$ is $\alpha$-stable. \qed
\end{prop}

\begin{rem}\label{sta-nil}
In the nilpotent case, the proposition remains true; here $\alpha=0$,
and, following the proof of lemma \ref{finite-deg}, the degree is
finite for all subsheaves with flat restriction to $T_\infty$.
\end{rem}

\begin{rem}
It is important to note that the stability condition just defined
is not an empty one. Indeed, $\alpha$-unstable bundles $\ee\seta\tproj$ 
can be obtained as extensions in the following way:
$$ 0 \seta p_1^*L_{\xi_0}\otimes p_2^*\oo_{\proj}(b) \seta \ee
   \seta p_1^*L_{-\xi_0}\otimes p_2^*\oo_{\proj}(-b)\otimes{\cal I}_k \seta 0 $$ 
where $b>0$ and ${\cal I}_k$ is the ideal sheaf of $k>0$ points in $\tproj$,
and we assume that none of these points are in $T_\infty$. Every sheaf
$\ee$ so obtained is locally-free, since the sheaf on the LHS is locally-free
and the one on the RHS is torsion-free. Clearly,
$\ee$ has trivial determinant, instanton number $k$ and asymptotic states $\as$.
\end{rem}

To finish the proof of theorem \ref{extn}, it remains to prove the
following proposition.

\begin{prop}\label{stable-implies-instanton}
Every $\alpha$-stable, holomorphic $S\ell_2$-bundle $\ee$ over $\tproj$ 
restricting to $L_{\xi_0}\oplus L_{-\xi_0}$ on $T_\infty$
can be obtained as the holomorphic extension of an instanton
on $T\times\cpx$ with asymptotic states $\as$, and whose monodromy 
around the torus at infinity has eigenvalues $\exp(\pm 2\pi i \alpha)$.
\end{prop}

\begin{proof}
We will give two different ideas to prove the proposition, but we will
not give the proofs, because they follow essentially well known
arguments.

The first idea is direct construction: construct a Hermitian-Einstein
metric on $\ee|_{T\times\cpx}$ (so that the Chern connection is
anti-self-dual); for this, one has first to build a metric $h_0$ on
$\ee$ which gives asymptotically at infinity an instanton: this is
possible because $\alpha$ and the behavior of $\ee$ near infinity (see
(\ref{dev-ee})) give all the parameters at infinity of the instanton;
then one wants to deform $h_0$ to a solution $h$ of the
Hermitian-Einstein equation, mutually bounded with $h_0$; Simpson's
method \cite{S0} cannot be used, because $T\times\cpx$ has infinite
volume, but one can apply the method in \cite{B3}, using precise analysis at
infinity, which will be explained in the next section for the study of
the moduli space.

The second idea, giving a different proof, consists in using the Nahm
transform of instantons. Recall that our instantons are in
correspondence with Higgs bundles with singularities on the dual torus
$\hat{T}$, with a harmonic metric. Actually, the correspondence has a
purely holomorphic interpretation, and this is an occurrence of the
so-called {\em Fourier-Mukai transform}. Stability is ususally preserved 
by such a correspondence, so that an $\alpha$-stable bundle on $\tproj$ 
would transform into a stable parabolic Higgs bundle on $\hat{T}$; then 
one can apply Simpson's theorem \cite{S} to construct a harmonic metric, 
whose inverse Nahm transform provides an instanton with quadratic curvature 
decay, and by theorem \ref{asymp.par} this instanton has exactly the desired 
behavior at infinity.
\end{proof}

%---------------------------------------------------

\section{Moduli spaces} \label{sec.moduli}

We now proceed to the differential geometric construction of the
moduli space. The $L^2$ metric will then provide a hyperk\"ahler
structure on it.

We will restrict to the semisimple case; this choice simplifies the
construction, because theorem \ref{asymp.par} says that it is enough
to look at functional spaces with weights which are powers of $r$; the
analysis in the nilpotent case is possible, as in \cite{B3}, but
requires functional spaces with logarithmic weights.

Recall the model connection on the bundle $E$, trivialized near infinity:
\begin{eqnarray*}
A_0 = d &+& 
i\left(\begin{array}{cc} \lambda_1 & 0 \\ 0 & -\lambda_1 \end{array}
 \right) dx +
i\left(\begin{array}{cc} \lambda_2 & 0 \\ 0 & -\lambda_2 \end{array}
 \right) dy \\ &+&
i\left(\begin{array}{cc} \mu_1\cos\theta-\mu_2\sin\theta & 0 \\
 0 & -\mu_1\cos\theta+\mu_2\sin\theta \end{array}
 \right) \frac{dx}{r} \\ &+&
i\left(\begin{array}{cc} \mu_1\sin\theta+\mu_2\cos\theta & 0 \\
 0 & -\mu_1\sin\theta-\mu_2\cos\theta \end{array}
 \right) \frac{dy}{r} + \\ &+&
i\left(\begin{array}{cc} \alpha & 0 \\ 0 & -\alpha \end{array}
 \right) d\theta.
\end{eqnarray*}
Note that in order to get $L^2$ deformations, we cannot move the
parameters $\lambda$, $\mu$ and $\alpha$; in view of theorem
\ref{asymp.par}, it is natural to consider connections $A_0+a$, such that
$$
|a|=O(r^{-(1+\delta)}), \quad |\nabla_{A_0}a|=O(r^{-(2+\delta)});
$$
actually, this $C^1$ space is not good for analysis, and we have the
choice to substitute either a H\"older space $C^{1,\eta}$ or a Sobolev
space $L^{1,p}$; we make the last choice, for $p$ big enough, and this
leads to the technical definitions
\begin{eqnarray*}
\Omega^1_\delta &=& \{ a\in\Omega^1(\mathfrak{su}(E)), a\in L^p_{1-2/p+\delta},
                          \nabla_{A_0}a\in L^p_{2-2/p+\delta} \} \\
\cA &=& A_0 + \Omega^1_\delta \\
\cG &=& \{ g\in SU(E), \nabla_{A_0}g g^{-1}\in\Omega^1_\delta\} \\
\cF &=& \{ F\in \Omega^2_+(\mathfrak{su}(E)), F\in L^p_{2-2/p+\delta} \}.
\end{eqnarray*}
The Lie algebra of $\cG$ is
$$ T_1\cG =\{ u\in\mathfrak{su}(E), \nabla_{A_0}u \in\Omega^1_\delta\}.$$
Note that for $a\in\Omega^1_\delta$, lemma \ref{l2} implies that
actually $a_\perp\in L^p_{2-2/p+\delta}$, so that this Sobolev space
is the same as the one considered is part \ref{part1}.
Also, the Sobolev embedding (\ref{Sobolev-embedding}) implies
$\Omega^1_\delta \subset C^0_\delta$, and an important property is
that the embedding $\Omega^1_{\delta} \subset C^0_{\delta'}$ is
compact if $\delta'<\delta$; gauge transformations $g\in\cG$ can be
continuously extended over $T_\infty$, so that
$$
g|_{T_\infty} = \left(\begin{array}{cc} u & 0 \\ 0 & u^{-1} \end{array}
                \right),
$$
where $u\in S^1$ is fixed.
Also, $\cG$ acts smoothly on $\cA$ and the 
curvature is a smooth map from $\cA$ to $\cF$.

Remark that there is no reducible connection in $\cA$, since a reduction would
decompose the bundle $E$ as $L\oplus L^{-1}$, with $L$ topologically
trivial on the torus at infinity; but then we would get $c_2(E)=0$.

Now we need the following proposition;
the proof is given at the end of the section.
\begin{prop}\label{Delta} For $k>0$ and $A\in\cA$, we have:
\begin{enumerate}
\item the laplacian $\Delta_A:T_1\cG\rightarrow
L^p_{2-2/p+\delta}$  is an isomorphism; therefore there is a slice at
$A$ to the action of $\cG$ on $\cA$, given by $\{A+a,d_A^*a=0\}$;
\item if $A$ is an instanton, then the map
$d_A^+\oplus d_A^*:\Omega^1_\delta\rightarrow L^p_{2-2/p+\delta}$ is Fredholm
surjective; the kernel coincide with the $L^2$-kernel.
\end{enumerate}
\end{prop}
Note that in the first statement of the proposition,
it was crucial to allow gauge transformations to take non
trivial values on $T_\infty$, otherwise one cannot obtain the slice
$d_A^*a=0$.

Define the moduli space $\cM$ as the space of instantons
$A\in\cA$ modulo the gauge group $\cG$.
As is well-known, $F_A^+$ is a hyperk\"ahler moment map for the action
of $\cG$ on $\cA$ with respect to the three complex structures on
$\torus$:
\begin{eqnarray}
I_1(z_1,z_2,w_1,w_2) & = & (-z_2,z_1,-w_2,w_1) \nonumber \\
I_2(z_1,z_2,w_1,w_2) & = & (-w_1,w_2,z_1,-z_2) \label{cpx.str1} \\
I_3(z_1,z_2,w_1,w_2) & = & (-w_2,-w_1,z_2,z_1) \nonumber
\end{eqnarray}
where $z=z_1+iz_2$ and $w=w_1+iw_2$. With the help of the previous proposition, 
standard theory now gives us:

\begin{prop}\label{moduli}
The moduli space $\cM$ is a smooth hyperk\"ahler manifold; the
tangent space at $[A]$ is isomorphic to the $L^2$-kernel of
$d_A^+\oplus d_A^*$ acting on $\Omega^1(\mathfrak{su}(E))$.
It has dimension $8k-4$.
\end{prop}

\begin{proof}[Proof of proposition \ref{Delta}]
First, we have to understand the behavior of the laplacian $\Delta_A$ acting
on sections of $End(E)$. We want to prove that it is Fredholm.
This property is not changed by a perturbation in $\Omega^1_\delta$
(this adds to $\Delta_A$ a compact operator),
and we can therefore restrict to the case when $A=A_0$ on $r\geq R$.
On this domain $r\geq R$, the laplacian preserves the decomposition
$u_\Gamma\oplus u_\perp$.

The case of $u_\perp$ is easier: since we
have seen that $\|\nabla_{A_0}^2u\|_{L^p_{2-2/p+\delta}}$ controls
$\|u\|_{L^p_{2-2/p+\delta}}$, it follows that $A_0-\Gamma$, which is
$O(r^{-1})$, is small if $R$ is big enough; therefore
lemma \ref{neumann-Lp} proves that $\nabla_A$ is an isomorphism on
$r\geq R$ for the Neumann boundary condition (the same is true for
Dirichlet boundary condition).

The case of $u_\Gamma$ is more complicated, but can be reduced to
standard theory: recall that $u_\Gamma$ is torus invariant, so that
the operator now reduces to an operator on $\bR^2$;
the action of $\Delta_\Gamma$ on off-diagonal
coefficients (which exist only when $\Gamma$ is trivial) is by
$$ \frac{1}{r^2} \big( -(r\partial_r)^2 + (\partial_\theta^2\pm 2i\alpha)^2
                       +|\mu|^2 \big), $$
and the action on diagonal coefficients is the standard laplacian on
$\bR^2$ (that we obtain by making $\alpha=\mu=0$ in the previous formula);
now $r^2\Delta_A$ becomes the translation invariant laplacian
$$ -\partial_t^2 - (\partial_\theta^2\pm 2i\alpha)^2+|\mu|^2 $$
on the conformal cylinder $\bR_+\times S^1$, so that standard theory
\cite{L-M} now applies: such operator (say, with Dirichlet boundary
condition on $r=R$) is Fredholm for all weights,
except a discrete set of critical weights $\delta$ (they are
characterized by the existence at infinity of solutions of type
$\exp(-\delta t)t^k$); moreover, as the operator is self-adjoint, its index
is 0 at the weight 0 if it is noncritical, or $-1$ for small positive
weights if 0 is critical; in our situation, $u\in T_1\cG$ corresponds to
the decay $u\in L^p_{\delta-2/p}$, and this becomes exactly the weight
$\delta$ on the cylinder; there are two cases: if $\alpha$ or $\mu$ is
non zero (off-diagonal coefficients), then the weight 0 is not
critical, and the operator remains Fredholm for nearby $\delta$, with
index 0: actually is is an isomorphism, because it easy to verify that
is has no kernel; if $\alpha$ and $\mu$ are zero, then the laplacian
has index $-1$ for small weights $\delta>0$, so that it becomes an
isomorphism if we add the possibility to consider solutions $u$ of
$\Delta u=v$ with $u$ having some nonzero limit at infinity
(and this is exactly our definition of $\cG$). All these results
can also be checked by direct calculation, after
decomposing $u$ into Fourier series along each circle.

Finally, we deduce from these considerations that the laplacian
$\Delta_{A_0}$ is an isomorphism $T_1\cG\rightarrow L^p_{2-2/p+\delta}$
for the Dirichlet boundary conditions on $r\geq R$, and gluing this
isomorphism with a parametrix on the compact part, it follows
that $\Delta_A$ is Fredholm on $T\times\bR^2$.

In order to calculate the index, if $\Gamma$ is nontrivial, we have
seen that the index is not changed if we modify $A$ so that $A=\Gamma$
near infinity; the index of a self-adjoint operator on a compact
manifold is zero; by an excision principle, this has the consequence
that the index comes only from the contribution at infinity;
therefore, it is equal to the index of the operator $\Delta_\Gamma$
acting on the trivial bundle $\mathfrak{su}(\cpx^2)$; now this
operator is completely explicit: on the $u_\perp$ component, it is an
isomorphism, and on the $u_\Gamma$ component (that is, diagonal, torus
invariant, matrices), it is simply the standard laplacian in $\bR^2$,
and its index between the spaces that we have defined is again 0, with
1-dimensional kernel and cokernel equal to constant diagonal matrices.

If $\Gamma$ is trivial, we cannot reduce to the operator of flat
space, but we can reduce to $\Delta_{A_0'}$, with $A_0'$ the diagonal
connection
\begin{equation}\label{A0p}
A_0' = \chi(r) A_0 + (1-\chi(r)) d ,
\end{equation}
where $\chi(r)$ is a cutoff function which equals 1 for $r>R$ and 0
for $r<R-1$; then, as above, it is not difficult to prove that
$\Delta_{A_0}'$ is an isomorphism on non-diagonal components (and the
operator on the diagonal components is the same as above).

Finally, the operator $\Delta_A$ has no kernel in $T_1\cG$, since an
element in the kernel would decompose $A$, which is impossible.
This finishes the proof of the first part of the proposition.

If $A\in\cA$ is an instanton, observe that the operator $d_A^+ d_A^*$
acting on self-dual 2-forms, by the Weitzenb\"ock formula, equals the
laplacian $\nabla_A^*\nabla_A$; this means that the above results
remain true for $d_A^+ d_A^*$, and we deduce that the operator
$$
d_A^+ \oplus d_A^* : \Omega^1_\delta \longrightarrow L^p_{2-2/p+\delta}
$$
is surjective; its kernel equals the kernel of the laplacian
$2(d_A^+)^*d_A^+ + d_A d_A^*$; again one can prove (in particular
using lemma \ref{defn.op}) that this operator is Fredholm (for the
weight $\delta$); remark that the $L^2$ condition corresponds to a
critical weight (on diagonal components, where the operator is
asymptotically the standard laplacian of $\bR^2$), when $\Omega^1_\delta$
corresponds to a slightly greater weight; nevertheless, it remains
true that the $L^2$-kernel equals the kernel for slightly greater weights
(the possible new solutions in the kernel at the critical weight are
never $L^2$).
\end{proof}

\begin{proof}[Proof of proposition \ref{moduli}]
It remains only to calculate the dimension, which, by proposition
\ref{Delta}, is the index of the operator $d^+_A\oplus d_A^*$.
If the limit flat connection $\Gamma$ is non trivial,
this is simple to calculate by comparison to the same operator for $\Gamma$:
actually, by the excision principle,
$$
ind(d^+_A\oplus d_A^*) = ind(d_\Gamma^+\oplus d_\Gamma^*)+8k ;
$$
now for the flat connection $\Gamma$, the operator $d_\Gamma^+\oplus
d_\Gamma^*$ has no kernel (by the Weitzenb\"ock formula), but its
cokernel equals the  cokernel of the operator
$d_\Gamma^* d_\Gamma + d_\Gamma^+ d_\Gamma^*$ acting on
$\Omega^0(\mathfrak{su}(E))\oplus\Omega^2_+(\mathfrak{su}(E))=\bR^4\otimes\mathfrak{su}(E)$;
we have seen above that the cokernel of this operator on $\mathfrak{su}(E)$
is the $L^2$-orthogonal of constant, diagonal matrices. This proves
the formula for the index.

If $\Gamma$ is trivial, the same result holds, but one must compare
with the operator $d_{A_0'}^+\oplus d_{A_0'}^*$ defined in (\ref{A0p}).
\end{proof}

\paragraph{Fibration structure.}
It was shown in \cite{J3} that the moduli space of rank two holomorphic 
vector bundles over $\tproj$ with trivial determinant and instanton 
number $k$ contains an open set $\cM^*_k$ (corresponding  to the so-called 
\emph{regular bundles}) which has the structure of a fibration:
$$ \mathbb{T} \cdots \cM^*_k \seta \Sigma_k $$
The fibres are complex tori of complex dimension $2k-1$, and the base can 
be interpreted as the set of rational maps $\proj\seta\proj$ of degree $k$,
so that ${\rm dim}\Sigma_k = 2k+1$.

Fixing the splitting of $\ee$ at $T_\infty$, i.e. fixing the asymptotic state 
of the corresponding instanton connection $A$, amounts to fixing the value of 
these rational maps at $\infty\in\proj$. Moreover, as we will see in the 
next section, fixing the residue of $A$ amounts to fixing the first 
derivative at $\infty\in\proj$. 

Therefore, according to theorem \ref{extn}, we conclude that $\cM_{(k,\as,\mu)}$,
the moduli space of $SU_2$ doubly-periodic instantons with fixed instanton number 
$k$, asymptotic states $\as$ and residue $\mu$ with the complex structure induced 
from the complex structure $I_1$ on $T\times\bR^2$, is a fibration over $\Sigma_{(k,\as,\mu)}$,
the space of rational maps $f:\proj\seta\proj$ with fixed $f(w=\infty)$ and
$f'(w=\infty)$, with fibres given complex tori of dimension $2k-1$. 

Moreover, it is possible to show that the such fibres are lagrangian with 
respect to complex symplectic structure on $\cM_{(k,\as,\mu)}$ induced from 
the complex symplectic structure $\omega_{I_2}+I_1\omega_{I_3}$ on $\torus$ 
(see \cite{JM} for the proof of a similar result for elliptic K3 and abelian 
surfaces). 

\paragraph{An example: $\mathbf{k=1}$.}
We shall now give an explicit model for the moduli space of doubly-periodic
instantons with $k=1$; clearly, we also assume that $\xi_0\neq-\xi_0$ and 
$\mu\neq0$.

Our approach is based on the observations made above, that is, we shall 
study the set of rational maps $f:\proj\seta\proj$ of degree 1; in a 
neighbourhhod of $\infty\in\proj$, such maps can be written as follows:
$$ f(w) = \frac{w+b}{cw+d},\ \ {\rm where}\ w=0\ {\rm corresponds\ to}\ \infty\in\proj. $$
As we discussed above, we must still fix $f(0)$ and $f'(0)$. This means that
$b/d$ and $(d-cb)/d^2\neq 0$ are fixed. Thus, $\Sigma_{(1,\as,\mu)}=\cpx$, so that 
$\cM_{(k,\as,\mu)}$ is an elliptic fibration over $\cpx$.

Actually, one can say more: there is an action of $T\times\cpx$ on the
moduli space (by translations), so the moduli space is exactly
$T\times\cpx$, and the metric is flat.

\part{Nahm transform} \label{part3}

We now shift our attention to the Nahm transform of doubly-periodic
instanton connections \cite{J2}. Note that this
transform was defined in \cite{J2} only for instantons such that the
restriction of the underlying holomorphic bundle to a generic torus is
$L_\xi\oplus L_{-\xi}$ (this is what we called the semisimple case).
In this part, we shall restrict to this case.

Throughout this part, we assume familiarity with \cite{J2}, but
let us quickly recall how Nahm transform is defined. Given an
instanton $A$ on a $SU_2$-bundle $E$ on $T\times\bR^2$, one may twist
$A$ by a flat connection on $T$; these twists $A_\xi$ are
parameterized by $\xi\in\dual$. Now there is a coupled Dirac operator
$$ D_{A_\xi} : \Gamma(S^+\otimes E) \longrightarrow \Gamma(S^-\otimes E) $$
and one can show that the bundle of $L^2$-cokernels of $D_{A_\xi}$ is a rank
$k$ vector bundle $V$ over $\dual\setminus\{\pm\xi_0\}$; there is a natural
connection $B$ on $V$ obtained by projection, and one can define an
endomorphism $\Phi$ of $V$ by taking an element $\beta\in\ker D_{A_\xi}^*$
to the projection of $w\beta$ on this kernel; the pair $(B,\Phi)$ satisfies
Hitchin's equations on $\dual\setminus\{\pm\xi_0\}$.

From the holomorphic point of view, the picture is very clear: the
spinor bundle $S$ is identified $\Lambda^{0,*}$, so that the $L^2$-kernel of
$D^*_{A_\xi}$ is exactly the $L^2$-kernel of
$\overline{\partial}_{A_\xi}\oplus\overline{\partial}_{A_\xi}^*$ on
$\Omega^{0,1}\otimes E$. It can be proven that this $L^2$-kernel
coincides with $H^1(T\times\proj,\ee\otimes L_\xi)$, where $\ee$ is
the holomorphic extension of $A$ on $T\times\proj$; this provides a
holomorphic extension of $V$ on the whole $\dual$; this extension has
degree $-2$, as can be checked by Riemann-Roch theorem for families.
Moreover, there is a natural interpretation for the Higgs field:
one has the identification
\begin{equation}\label{decVw}
 V_\xi = H^1(T\times\proj,\ee\otimes L_\xi)
         = \bigoplus_{w\in\cpx} H^0(T_w,E\otimes L_\xi)
\end{equation}
where of course there is only a finite number of points $w\in\cpx$
(actually $k$, counted with multiplicity) such that $H^0(T_w,E\otimes
L_\xi)\neq 0$. Now the Higgs field $\Phi$ is multiplication by $w$ on 
$H^0(T_w,E\otimes L_\xi)$. From this description, one can see that the
Higgs field has a simple pole at $\pm\xi_0$ with semisimple residue,
and the residue has only one nonzero eigenvalue if $\xi_0\neq -\xi_0$,
two otherwise.

We first study how the new asymptotic parameters of doubly-periodic
instantons introduced in Part II behave under Nahm Transform. This
will prepare the way for the proof of theorem \ref{iso.thm}, our
last result.

\section{Asymptotic parameters} \label{ap.nt}

Following the general philosophy that the Nahm Transform is a sort
of nonlinear Fourier Transform, it is reasonable to expect the
asymptotic behavior of the instanton to be translated into further
singularity data for the Higgs field.

Recall that $\ee|_{T_\infty}=L_{\xi_0}\oplus L_{-\xi_0}$. From
(\ref{decVw}) we deduce a holomorphic splitting of $V$ on a small neighborhood
of $\as$:
\begin{equation}\label{decVBR}
V_\xi = B_\xi \oplus R_\xi
\end{equation}
where $B_\xi$ corresponds to the points in $\cpx$ that remain bounded
as $\xi\seta\xi_0$ and $R_\xi$ corresponds to the points that go off
to infinity. Clearly, $B_\xi$ approaches the kernel of ${\rm Res}_{\as}\Phi$
as $\xi\seta\xi_0$, while $R_\xi$ approaches the eigenspace
of the nontrivial eigenvalues of the residue.

The behaviour of the Higgs bundle with harmonic metric near the
singularities $\pm\xi_0$ is completely determined by the following theorem.

\begin{thm} \label{ap.nt.thm}
Let $A$ be a doubly-periodic instanton with limiting holonomy
$\alpha$ and residue $\mu$; let $(B,\Phi)$ be its Nahm transformed
Higgs pair.
The unique nonzero eigenvalue of ${\rm Res}_{\pm\xi_0}\Phi$ is given
by $\pm\mu$.
In the decomposition (\ref{decVBR}), the harmonic metric on $V$
remains bounded on $B$, but behaves like $|\xi\pm\xi_0|^{1\pm\alpha}$
on $R$.
\end{thm}

\begin{rem} The sum of the degree of $V$, that is $-2$, and of the
weights $1\pm\alpha$, equals 0, as must be for a solution of Hitchin's
equations. The monodromy of the connection $B$ near the punctures is
semisimple, with only one nontrivial eigenvalue $\exp(\mp 2\pi i\alpha)$
on $R$ (or two if $\xi_0=-\xi_0$).
\end{rem}

We first prove the statement concerning the residues. The argument
to establish the statement concerning the limiting holonomy is
much more technical, and will involve a series of lemmas.

\paragraph{Residues.}
Let $\rho=r^{-1}$ and let $w'=w^{-1}=\rho e^{-i\theta}$ be a coordinate
near $\infty\in\proj$. Clearly, the holomorphic structure on the restriction
$E|_{T_{w'}}$ is given by the $(0,1)$-part of the $A|_{T_{w'}}$. Rewriting
equation (\ref{dev-ee}) in terms of $w'$, we obtain:
$$ \del_A |_{T_{w'}} = \del +
\left(\begin{array}{cc} \lambda & 0 \\ 0 & -\lambda \end{array}\right) d\overline{z}
+ \left(\begin{array}{cc} \mu & 0 \\ 0 & -\mu \end{array}\right) w'd\overline{z}
+ O(\rho^2) $$
so that:
$$ \left. \frac{d}{dw'} \big( \del_A|_{T_{w'}} \big) \right|_{w'=0} =
\left(\begin{array}{cc} \mu & 0 \\ 0 & -\mu \end{array}\right) $$
In other words, the residue $\mu$ can be regarded as the
infinitesimal variation of the holomorphic bundle $\ee|_{T_{w}}$
at $w=\infty$.

Since for every $w'$ sufficiently close to $\infty\in\proj$ we can
assume that \linebreak $E|_{T_{w'}}=L_{\xi(w')}\oplus L_{-\xi(w')}$, 
the above expression implies that:
$$ \left. \frac{d}{dw'} \xi(w') \right|_{w'=0} = \mu . $$
The eigenvalue of $\Phi$ going to infinity is $w(\xi)=1/w'(\xi)$ by
(\ref{decVw}); the statement follows.\qed

\paragraph{Limiting holonomy.}
Let us now look at the coupled Dirac laplacian $\Delta_{A_\xi}$ acting
on sections of $S^+\otimes E$; since $A$ is an instanton, we have that
\linebreak $D_{A_\xi}^*D_{A_\xi}=\nabla_{A_\xi}^*\nabla_{A_\xi}$, 
i. e. the Dirac laplacian coincides with the trace laplacian.
This laplacian is inversible in $L^2$ for $\xi\neq\pm\xi_0$ (see
\cite{J2}; this is also a consequence of the lemmas below), and we
note its inverse by $G_{A_\xi}$.
Such inverse is useful to produce harmonic representative of elements
of $H^1(T\times\proj,\ee\otimes L_\xi)$. Indeed, if we have a compactly
supported (0,1)-form $\beta$ with values in $E$ such that
$\overline{\partial}_{A_\xi}\beta=0$, then the $L^2$-harmonic representative
of the class $[\beta]$ is given by
$$ \beta - \overline{\partial}_{A_\xi} G_{A\_xi}
           \overline{\partial}_{A_\xi}^* \beta . $$

We now want to understand the inverse $G_{A_\xi}$ when $\xi$ approaches
the asymptotic states $\pm\xi_0$.
For simplicity, assume that $\xi_0=0$ in the next three lemmas; the
general case can be obtained by substituting $\xi-\xi_0$ for $\xi$
in the expressions below.

We know from theorem \ref{asymp.par}, where $\lambda=\xi$:
$$ A_\xi = (A_0)_\xi + a \ \ \ {\rm with}\ \ |a| = O(r^{-1-\epsilon}) $$
and
\begin{eqnarray*}
(A_0)_\xi & = & d+i\left( \begin{array}{cc} \alpha & 0 \\ 0 & -\alpha \end{array} \right)d\theta +
i\left( \begin{array}{cc} \lambda_1dx+\lambda_2dy & 0 \\
        0 & -\lambda_1dx-\lambda_2dy \end{array} \right) + \\
& & \phantom{d+} \frac{i}{r}\left( \begin{array}{cc} \mu_1dx+\mu_2dy & 0 \\
                  0 & -\mu_1dx-\mu_2dy \end{array} \right).
\end{eqnarray*}
We assume also that at either $\mu_1$ or $\mu_2$ is nonzero;
however, the proofs below will also work if $\mu_1=\mu_2=0$, but
$\alpha\neq0$.

\begin{lem} \label{l1p3}
Let $\sigma$ is a section of $E\seta\torus$; if $\lambda$ is
sufficiently small and $|w|$ is large enough, then:
$$ \int_{T_w} |\nabla_{(A_0)_\xi}\sigma|^2 \geq
   \left| \lambda + \frac{\mu}{w} \right|^2 \int_{T_w} |\sigma|^2 .$$
\end{lem}

\begin{proof}
Consider the Fourier expansion $\sigma=\Sigma \sigma_{nm} e^{i(nx+my)}$.
Then on the torus $T_w$, we have:
\begin{eqnarray*}
\int_{T_w} |\nabla\sigma|^2 &=&
\int_{T_w}
\left| \big(\partial_x+i\lambda_1+i\frac{\mu_1}{|w|}\big)\sigma \right|^2+
\left| \big(\partial_y+i\lambda_2+i\frac{\mu_2}{|w|}\big)\sigma \right|^2\\
& = & \sum \left| (n+im+\lambda+\frac{\mu}{|w|} \right|^2 |\sigma_{nm}|^2 .
\end{eqnarray*}
However, under the hypothesis above,
$$ \left| n+im+\lambda+\frac{\mu}{|w|} \right| \geq
   \left| \lambda + \frac{\mu}{w} \right| $$
for all $n,m$, which proves the lemma.
\end{proof}

\begin{lem} \label{l2p3}
Under the hypothesis of lemma \ref{l1p3}, we have:
\begin{equation} \label{p3est1}
\int_{r\geq R} |\nabla_{(A_0)_\xi}\sigma|^2 \geq
   c |\mu|^2 \int_{r\geq R} \frac{|\sigma|^2}{r^2}
\end{equation} \end{lem}

\begin{proof}
By the previous lemma, the estimate holds away from the region
where $|\lambda+\mu/w|$ is small, that is:
$$ \frac{1}{2}\frac{|\lambda|}{|\mu|} \leq |w| \leq 2\frac{|\lambda|}{|\mu|} $$
Actually, we claim that if the estimate of the lemma is satisfied outside this
region, then it must be satisfied everywhere. Indeed, one has the
inequality for any function $f:\real^2\seta\real$, and a constant $c$
independent of $\rho$,
\begin{equation} \label{ineq}
\int_{\rho\leq r\leq 2\rho} \frac{f^2}{r^2} \leq c\cdot
\left( \int_{2\rho\leq r\leq 4\rho} \frac{f^2}{r^2} +
\int_{\rho\leq r\leq 4\rho} |\partial_r f|^2 \right)
\end{equation}
and the lemma follows by applying (\ref{ineq}) to $f=|\sigma|$ and
$\rho=|\lambda/\mu|$. The proof of (\ref{ineq}) is left to the
reader.
\end{proof}

Note that an estimate similar to (\ref{p3est1}) remains valid if $\mu=0$,
but $\alpha\neq0$. In fact, the proof is even simpler, since one
has the estimate:
$$ \int_{r=R} |\nabla_{(A_0)_\xi}\sigma|^2 \geq
   \frac{|\alpha|^2}{r^2} \int_{r\geq R} |\sigma|^2 $$
from which one immediately obtains:
\begin{equation} \label{p3est2}
\int_{r\geq R} |\nabla_{(A_0)_\xi}\sigma|^2 \geq
   |\alpha|^2 \int_{r\geq R} \frac{|\sigma|^2}{r^2} .
\end{equation}

\begin{lem} \label{l3p3}
The solution of the Poisson equation $\Delta_{A_\xi}u=v$ satisfies:
$$ \| r^{-1}  u \|_{L^2} + \|\nabla_{A_\xi}u \|_{L^2} \leq c \|r  v \|_{L^2} $$
$$ {\rm and}\ \ \
   |\xi|^2 \| u \|_{L^2} + |\xi|\cdot\|\nabla_{A_\xi}u \|_{L^2} \leq c \|r  v \|_{L^2} $$
\end{lem}

\begin{proof}
First, note that:
\begin{equation} \label{p3est3}
\int |\nabla_{A_\xi}\sigma |^2 \geq c
\left( |\xi|^2 \int |\sigma|^2 + \int \frac{|\sigma|^2}{r^2} \right) .
\end{equation}
Near infinity, this a consequence of lemma \ref{l2p3} and of the fact that
$A=A_0+O(r^{-1-\epsilon})$. Globally, the estimate follows from
the Poincar\'e-type inequality:
$$ \int_{r\leq R} |\sigma|^2 \leq c
\left( \int_{r\leq R} |\nabla\sigma|^2 + \int_{R/2\leq r\leq R} |\sigma|^2 \right)$$

To prove the lemma itself, we have that:
\begin{eqnarray*}
\|\nabla_{A_\xi}u \|^2_{L^2} & = & \int \langle \Delta_{A_\xi}u,u \rangle =
\int \langle v , u \rangle \leq \\
& \leq & \| r  v \|_{L^2} \| r^{-1}  u \|_{L^2}
         \leq c \| r  v \|_{L^2} \| \nabla_{A_\xi}u \|_{L^2}
\end{eqnarray*}
by (\ref{p3est3}). Thus, we conclude that
$\| \nabla_{A_\xi}u \|_{L^2}\leq c\| r  v \|_{L^2}$, and again
by (\ref{p3est3}) we have $\| r^{-1}  u \|_{L^2}\leq c\| r  v \|_{L^2}$.
The second estimate is obtained in a similar way.
\end{proof}

We are now finally ready to complete the proof of theorem \ref{ap.nt.thm}.
Let us first analyze the behavior of the harmonic metric on the
local sub-bundle $B\hookrightarrow V$ with fibers given by $B_\xi$. Let
$\beta$ be a section of $B$. Then, for each $\xi\neq\xi_0$, we
know from (\ref{decVw}) that $\beta(\xi)$ can be represented as a
section of $\Lambda^{0,1}E\otimes L_\xi$ supported on $r\leq R$
for some $R$ sufficiently large. Furthermore, its harmonic
representative in $H^1(\torus,E\otimes L_\xi)$ is given by
$\beta(\xi)-\del_{A_\xi}G_{A_\xi}\del^*_{A_\xi}\beta(\xi)$.
By lemma \ref{l3p3}, we have:
$$ \| \del_{A_\xi}G_{A_\xi}\del^*_{A_\xi}\beta(\xi) \|_{L^2} \leq
   c \|r \del_A^*\beta(\xi) \|_{L^2} \leq c  R \| \del_A^*\beta(\xi) \|_{L^2} $$
which remains bounded even as $\xi\seta\xi_0$. This means that
the limit
$$ \beta(\xi_0)=\lim_{\xi\seta\xi_0} \beta(\xi) $$
has a square-integrable harmonic representative, so that the
harmonic metric restricted to the sub-bundle $B$ extends across
$\as$.

Now let $R\hookrightarrow V$ be a local sub-bundle with fibers
given by $R_\xi$; remind that near infinity, we have
$\ee|_{T_w}=L_{\xi(w)}\oplus L_{-\xi(w)}$; take a section $\beta(\xi)$
of $R_\xi$ coming by (\ref{decVw}) from sections of
$\ee|_{T_{w(\xi)}}\otimes L_\xi$ converging to a section of
$\ee|_{T_\infty}\otimes L_{\xi_0}=L_{2\xi_0}\oplus\cpx$.
Here we have to be more specific: say that a section $\sigma\in
H^0(T_{w(\xi)},\ee\otimes L_\xi)$ corresponds to the class in
$H^1(T\times\proj,\ee\otimes L_\xi)$ represented by the (0,1)-current
\begin{equation}\label{current}
\sigma(z) \delta_{w(\xi)}(w) d\overline{w} ,
\end{equation}
where $\delta_{w(\xi)}$ is the Dirac function at the point $w(\xi)$.
From this description, we see that,
for each $\xi\neq\xi_0$, the representative $\beta(\xi)$ can
be chosen with compact support near $r=|w(\xi)|$, and bounded in $L^{1,2}$.
Now lemma \ref{l3p3} gives, as above,
$$ \| \del_{A_\xi}G_{A_\xi}\del^*_{A_\xi}\beta(\xi) \|_{L^2} \leq
   c \|r \del_A^*\beta(\xi) \|_{L^2} \leq
   \frac{c}{|\xi-\xi_0|} \| \del_A^*\beta(\xi) \|_{L^2} . $$
This means that the norm of the harmonic representative of $\beta(\xi)$ is
bounded by $|\xi-\xi_0|^{-1}$.

This result must be interpreted, since (\ref{current}) actually does
not extend to $w=\infty$, so that our $[\beta(\xi)]$ is not a section
of $R$ which extends over the puncture $\xi_0$. There are two changes
to make; first, note that a (0,1)-form smooth on $\proj$
near infinity is $d\overline{w}/\overline{w}^2$, so we see that we must
consider $\beta(\xi)/\overline{w}(\xi)^2$ instead of $\beta(\xi)$. The second
change to be made is that we want $\beta(\xi)$ holomorphic in $\xi$.
This involves a constraint on the choice of $\sigma$: from the growth
of the holomorphic sections of $\ee$ at infinity studied in section
\ref{holo.extn}, it follows that $|\sigma|\sim |w(\xi)|^\alpha$, and we can
finally conclude that the norm of a holomorphic section of $R$ is bounded by
$|\xi-\xi_0|^{1-\alpha}$.

From these results, it follows that the harmonic metric of the Higgs
bundle $V$ extends on $B$, and is bounded by
$|\xi\pm\xi_0|^{1\pm\alpha}$ on $R$. This gives a bound $1\pm\alpha$
for the weights of the parabolic structure of $V$.
However, the ``parabolic degree'' of the bundle must be zero, and $V$
has degree $-2$, so that the weights must be exactly equal to $1\pm\alpha$.
\qed

\paragraph{Reformulating the Nahm transform theorem}
Together with \cite{J1,J2}, theorem \ref{ap.nt.thm} allows us to
state a complete version of the Nahm transform theorem, including
the new asymptotic parameters defined in Part II:

\begin{thm} \label{nahm}
The Nahm transform is a correspondence between the following
objects:
\begin{itemize}
\item $SU(2)$ doubly-periodic instantons with instanton
number $k>0$ and asymptotic parameters $(\as,\alpha,\mu)$;
\item rank $k$ logarithmic Higgs bundles with harmonic metric over $\dual$
with singularity behavior as described in theorem \ref{ap.nt.thm}.
\end{itemize} \end{thm}

%---------------------------------------------------------------

\section{The hyperk\"ahler property} \label{iso}

Our final task is to prove that the Nahm transform of doubly-periodic
instantons define a hyperk\"ahler isometry between $\cM$, the
moduli space of doubly-periodic instanton constructed in section
\ref{sec.moduli}, and $\hmod$, the moduli space of meromorphic Higgs
pairs satisfying the conditions of theorem \ref{nahm}. To do that,
we shall follow the  following strategy. First, we compute the
derivative of the map:
\begin{eqnarray*}
N:\ \cM & \longrightarrow & \hmod \\
A & \mapsto & (B,\Phi)
\end{eqnarray*}
defined by the Nahm transform, verifying that it is indeed well-defined.
We then show that $D_{[A]}N$ preserves the three complex structures in
each space. The last step is to show that $D_{[A]}N$ preserve the
metrics in each space.

\paragraph{Computing the derivative.}
Recall the definition of the tangent space $T_{[A]}\cM$ at the gauge
equivalence class of an instanton $A$ can be characterized as follows:
\begin{equation}
T_{[A]}\cM = \left\{ a\in L^2(\Omega^1\mathfrak{su}(E))\ {\rm s.t.}\ 
\begin{array}{rl} (i) & d_A^* a = 0 \\ (ii) & d_A^+ a = 0 \end{array} \right\}
\end{equation}
The 1-form $a$ is regarded as a infinitesimal variation of the
instanton connection $A$, inducing a 1-parameter family of
connections $A_t=A + t  a$, which are anti-self-dual up to
first order.

Now let $\{\Psi(\xi)^j\}_{j=1}^k$ be an orthonormal base for
coupled adjoint Dirac operator ${\rm ker}D_{A_\xi}^*$. In order to
compute the derivative $D_{[A]}N$, we must understand the
infinitesimal change on harmonic spinors induced by the
infinitesimal change on the instanton. We are looking for negative
spinors $\varphi(\xi)^j$ such that the 1-parameter family
$\Psi_t(\xi)^j=\Psi(\xi)^j+t\cdot\varphi(\xi)^j$ satisfies
$D_{(A_\xi)_t}^*\Psi_t(\xi)^j=0$ up to first order.
In other words,
$$ \left. \frac{d}{dt} D_{(A_\xi)_t}^*\Psi_t(\xi)^j \right|_{t=0} =
   D_{A_\xi}^*\varphi(\xi)^j + a\bullet\Psi(\xi)^j = 0 $$
where $\bullet$ means Clifford multiplication. Therefore, the
infinitesimal variations on harmonic spinors are given by:
\begin{equation} \label{infvar.spinors}
\varphi(\xi)^j = - D_{A_\xi}G_{A_\xi}(a\bullet\Psi(\xi)^j)
\end{equation}

Recall from \cite{J2} that the Nahm transformed Higgs pair is
defined as follows:
\begin{equation} \label{bphi} \begin{array}{ccc}
B(\xi)^{ij} = \langle \Psi(\xi)^i, \hat{d}\Psi(\xi)^j \rangle
& \ {\rm and}\ &
\Phi(\xi)^{ij} = \langle \Psi(\xi)^i, w\Psi(\xi)^j \rangle d\xi
\end{array} \end{equation}
where $\hat{d}$ means differentiation with respect to $\xi$, the
coordinate on the dual torus $\dual$, and the inner products are
taken in $L^2(E\otimes S^-)$. Thus, the infinitesimal change
in the Nahm transformed Higgs pair $(B,\Phi)$ is given by:
\begin{eqnarray}
b(\xi)^{ij} & = &
\left. \frac{d}{dt} \langle \Psi_t(\xi)^j, \hat{d}\Psi_t(\xi)^j \rangle \right|_{t=0} =
\nonumber  \\ & & \nonumber \\
& = & \langle G_{A_\xi}\Psi(\xi)^i, \Omega\bullet a\bullet \Psi(\xi)^j \rangle -
      \langle \Omega\bullet a\bullet \Psi(\xi)^i, G_{A_\xi}\Psi(\xi)^j \rangle
      \label{Bdefn}
\end{eqnarray}
and
\begin{eqnarray}
\phi(\xi)^{ij} & = &
\left. \frac{d}{dt} \langle \Psi_t(\xi)^j, w\Psi_t(\xi)^j \rangle \right|_{t=0} =
\nonumber \\ & & \nonumber \\
& = & \langle G_{A_\xi}\Psi(\xi)^i, dw\bullet a\bullet \Psi(\xi)^j \rangle d\xi
\label{Phidefn}
\end{eqnarray}
where $\Omega=i\big( d\xi_1dz_1 + d\xi_2dz_2 \big)$ is the curvature of the Poincar\'e bundle
over $T\times\dual$.

The tangent space $T_{[(B,\Phi)]}\hmod$ at the gauge equivalence class of a
Higgs pair $(B,\Phi)$, can described as follows (see for instance \cite{H}):
\begin{equation} \label{higgstgt}
T_{[(B,\Phi)]}\hmod = \left\{
\begin{array}{l}
b\in L^2(\Omega^1\mathfrak{u}(V)) \\
\phi\in L^2(\Omega^{1,0}\mathfrak{gl}(V)) \end{array}\ {\rm s.t.}\
\begin{array}{rl}
(i) &   d_B b + [\Phi,\phi^*] + [\phi,\Phi^*] = 0 \\
(ii) &  \del_B \phi + [b^{0,1},\Phi] = 0 \\
(iii) & d_B^* b + {\rm Re}[\Phi^*,\phi] = 0
\end{array} \right\}
\end{equation}
Again, $(b,\phi)$ define a 1-parameter family of pairs
$(B_t = B + t  b,\Phi_t = \Phi + t \phi)$ which satisfy
Hitchin's equations up to first order.

Therefore, it is clear from (\ref{Bdefn}) and (\ref{Phidefn}) that
the pair $(b,\phi)$ satisfies the linearized Hitchin's equations
($(i)$ and $(ii)$ in (\ref{higgstgt})).

We must only verify that $(b,\phi)$ are transversal to infinitesimal
changes in $(B,\Phi)$ arising from infinitesimal gauge transformations,
i.e. must check equation $(iii)$ in (\ref{higgstgt}). To do that,
denote by $\tilde{B}$ and $\tilde{b}$ the $(\real^2)^*$-invariant
1-forms on $\dual\times(\real^2)^*$ obtained from $(B,\Phi)$ and
$(b,\phi)$, respectively. Clearly, $\tilde{B}$ is anti-self-dual and
$$ d_B^* b + {\rm Re}[\Phi^*,\phi] = 0 \Leftrightarrow
   d_{\tilde{B}}^*\tilde{b}=0 $$
The following result completes our first step towards the proof
of theorem \ref{iso.thm}

\begin{lem}
If $d_A^* a = 0$, then $d_{\tilde{B}}^*\tilde{b}=0$.
\end{lem}
\begin{proof}
See proposition 3.1 in \cite{BVB}.
\end{proof}

\begin{rem}
Using the ideas above, one can easily compute the derivative of
the inverse Nahm transform, thus showing that $N:\cM\seta\hmod$
is a diffeomorphism. Noting that, since$\cM$ is smooth, the
diffeomorphism type of the moduli space of instantons does not
depend on the choice of asymptotic parameters $(\alpha,\lambda,\mu)$,
one concludes that the diffeomorphism type of the moduli of Higgs
bundles is independent not only of the singularity data (residues
and parabolic structure), as it was observed by Nakajima in \cite{N2},
but also of the position of the singularities.
\end{rem}

\paragraph{Commuting with the complex structures.}
Consider coordinates \linebreak $(\xi_1,\xi_2,\omega_1,\omega_2)$ on $(\real^4)^*$,
which are dual to $(z_1,z_2,w_1,w_2)$. Each of the complex structures (\ref{cpx.str1})
in $\real^4$ naturally induces a similar complex structures $\hI_j$ on
$(\real^4)^*$. Thus, we have maps:
$$ \Lambda^1\real^4 \stackrel{I_j}{\seta} \Lambda^1\real^4
\ \ \ {\rm and} \ \ \
\Lambda^1(\real^4)^* \stackrel{\hI_j}{\seta} \Lambda^1(\real^4)^* $$

The complex structures on $\hmod$ can be then defined as follows. As above,
let $\tilde{b}\Lambda^1(\real^4)^*$ be the $(\zed^2\times\real^2)^*$-invariant
1-form obtained from $(b,\phi)$. Then  $\hI_j(\tilde{b})$ is also a
$(\zed^2\times\real^2)^*$-invariant 1-form on $(\real^4)^*$, which can then be
interpreted as an element of (\ref{higgstgt}). It is easy to see that these
coincide with the complex structures originally defined by Hitchin in \cite{H}.
Therefore, we have to show that the following diagram:
\begin{equation} \label{diag1} \xymatrix{
\Lambda^1\real^4\su_2 \ar[r]^-{D_{[A]}N} \ar[d]^-{I_j} & 
\Lambda^1(\real^4)^*\ou_k \ar[d]^ -{\hI_j} \\
\Lambda^1\real^4\su_2 \ar[r]^-{D_{[A]}N} & \Lambda^1(\real^4)^*\ou_k }
\end{equation}
commutes. The horizontal maps are defined as follows:
\begin{equation} \label{btilde}
D_{[A]}N(a) = \tilde{b} =
\langle G_{A_\xi}\Psi(\xi)^i, \widetilde{\Omega}\bullet a\bullet \Psi(\xi)^j \rangle -
\langle \widetilde{\Omega}\bullet a\bullet \Psi(\xi)^i, G_{A_\xi}\Psi(\xi)^j \rangle
\end{equation}
with $\widetilde{\Omega}=i\big( d\xi_1dz_1 + d\xi_2dz_2 + d\omega_1dw_1 + d\omega_2dw_2 \big)$.

Each $I_j$ induces an isomorphism $l_j:\real^4\seta\cpx^2$ satisfying the following
commutative diagram:
\begin{equation} \label{diag2} \xymatrix{
\Lambda^1\real^4\su_2 \ar[r]^-{l_j} \ar[d]^-{I_j} & 
\Lambda^{(1,0)}\cpx^2\osl_k \ar[d]^-{\cdot i} \\
\Lambda^1\real^4\su_2 \ar[r]^-{l_j} & \Lambda^{(1,0)}\cpx^2\osl_k }
\end{equation}
where the map on the left hand side is multiplication by $i=\sqrt{-1}$. Of course,
a similar diagram holds for $\hat{l_j}:(\real^4)^*\seta(\cpx^2)^*$.

The key point is to note that each map:
$$ D_{[A]}N_{\cpx} = \hat{l_j} \circ D_{[A]}N \circ l_j^{-1}:
\Lambda^{(1,0)}\cpx^2 \seta \Lambda^{(1,0)}(\cpx^2)^* $$
$$ D_{[A]}N_{\cpx}(\alpha) = 
\langle G_{A_\xi}\Psi(\xi)^i, \widetilde{\Omega}_{\cpx}\bullet\overline{\alpha}\bullet \Psi(\xi)^j \rangle -
\langle \widetilde{\Omega}_{\cpx}\bullet\alpha\bullet \Psi(\xi)^i, G_{A_\xi}\Psi(\xi)^j \rangle $$
is $\cpx$-linear, where $\widetilde{\Omega}_{\cpx}=\hat{l_j}\times l_j(\widetilde{\Omega})$.
Therefore, we conclude:
\begin{eqnarray*}
\hI_j(D_{[A]}N(a)) & = & \hat{l_j}^{-1} \circ (\cdot i) \circ \hat{l_j} \circ D_{[A]}N(a) =
D_{[A]}N \circ l_j^{-1} \circ (\cdot i) \circ l_j (a) = \\
& = & D_{[A]}N (I_j(a))
\end{eqnarray*}
as desired.

\paragraph{The Nahm transform is an isometry.} 
Again, the fact that the Nahm transform is an isometry is actually a 
property of the underlying four-dimensional transform. The calculations
of Braam and van Baal \cite{BVB} are quite precise and also apply to 
the present situation.

Recall that the metric on the instanton moduli space is given by the
$L^2$ norm of the tangent vectors, that is:
$$ g(a_1,a_2) = \int_{\torus} {\rm Tr}(a_1\wedge*a_2) $$
while the metric on the Higgs moduli space is given by
$$ \hat{g}\big((b_1,\phi_1),(b_2,\phi_2)\big) =
   \int_{\dual} {\rm Tr}(b_1^*b_2 + \phi_1\phi_2^*) $$
or, equivalently, in terms of the 4-dimensional 1-forms
$\tilde{b_1}$ and $\tilde{b_2}$:
$$ \hat{g}(\tilde{b_1},\tilde{b_2}) = 
   \int_{\real^2}^* {\rm Tr}(\tilde{b_1}\wedge*\tilde{b_2}) $$
where integration is now done only with respect to the 
two coordinates on $(\real^4)^*$ on which $\tilde{b_1}$ and 
$\tilde{b_1}$ depend.

Let $(b,\phi)=D_{[A]}N(a)$; it is enough to show that:
$$ \hat{g} \big( D_{[A]}N(a) , (b,\phi) \big) =
   g \big( a , D_{[A]}N^{-1}(b,\phi) \big) $$
This can be done exactly as proposition 3.2 of \cite{BVB}.

Alternatively, we can reduce the isometry property to a purely algebraic 
statement as follows.

Fix the complex structure $I_1$ on $T^2\times\real^2$. The instanton moduli 
space $\cM$ is then identified with the moduli space of $\alpha$-stable 
holomorphic vector bundles $\ee\seta\tproj$ as a K\"ahler manifold. Moreover, 
its tangent space becomes identified with $H^1(\tproj,{\rm End}\ee)$. One can 
define a complex symplectic structure on $\cM$ via the bilinear pairing:
$$ H^1(\tproj,{\rm End}\ee) \times H^1(\tproj,{\rm End}\ee) 
   \stackrel{\omega}{\seta} H^2(\tproj,{\rm End}\ee) = \cpx $$

On the other hand, the moduli space of Higgs pairs $\hmod$ becomes identified,
as a K\"ahler manifold, with the moduli space of stable parabolic Higgs bundles. 
The tangent is then given by the hypercohomology $\hh^1$ of the following complex 
of sheaves: 
$$ ParEnd(\vv) \stackrel{[\cdot,\phi]}{\seta} \Lambda^1\otimes ParEnd(\vv) $$
where $ParEnd(\vv)$ is the sheaf of parabolic endomorphism of the holomorphic 
Higgs bundle $\vv$, see \cite{B3} for a detailed explanation. A complex symplectic 
structure on $\hmod$ can be defined via the bilinear pairing
$$ \hh^1 \times \hh^1 
   \stackrel{\hat{\omega}}{\seta} \hh^2 = \cpx $$

In order to show that the Nahm transform is an isometry, it is
enough to prove that the holomorphic version of the Nahm transform
(see \cite{J3}) preserves the bilinear pairings above. This is an 
algebraic statement, which one can hope to prove using spectral 
sequences. Indeed, as we mentioned before, the holomorphic version of
the Nahm transform of doubly-periodic instantons is an example of a 
{\em Fourier-Mukai transform}, which usually preserves this type of
pairings.

 \end{document}